\documentclass[10 pt, journal]{ieeetran}
\IEEEoverridecommandlockouts                              % This command is only needed if 
\usepackage{mydef}
\usepackage{cases}
\title{\LARGE \bf
A Proximal Diffusion Strategy for Multi-Agent Optimization with Sparse Affine Constraints
}

\author{Sulaiman A. Alghunaim$^{*}$, Kun Yuan$^{*}$,  and Ali H.~Sayed$^{\dagger}$, {\em Fellow}, IEEE% <-this % stops a space
\thanks{A short preliminary conference version appears in \cite{alghunaim2018cdc}. No convergence proofs were included in \cite{alghunaim2018cdc}. Besides proofs and derivations, this extended version also deals with the case of non-differentiable regularizres.}
\thanks{$^{*}$S. A. Alghunaim and K. Yuan are with the Electrical and Computer Engineering Department, University of California at Los Angeles (UCLA), CA 90095. emails:{\tt \{\small salghunaim,kunyuan\}@ucla.edu}.
        }%
\thanks{$^{\dagger}$A. H. Sayed is with the Ecole Polytechnique Federale de Lausanne
EPFL, School of Engineering, CH-1015 Lausanne, Switzerland email:
        {\tt\small ali.sayed@epfl.ch}. This work  was  supported  in  part  by  NSF  grants  CCF-1524250.}%
}

\begin{document}

%\address{$^{\star}$Department of Electrical and Computer Engineering,  University of California, Los Angeles \\
%$^{+}$School of Engineering, Ecole Polytechnique Federale de Lausanne, Switzerland}

\maketitle
\thispagestyle{plain}
\pagestyle{plain}

\begin{abstract}
This work develops a proximal primal-dual decentralized strategy for multi-agent optimization problems that involve multiple coupled affine constraints, where each constraint may involve only a subset of the agents. The constraints are generally sparse, meaning that only a small subset of the agents are involved in them. This scenario arises in many applications including decentralized control formulations, resource allocation problems, and smart grids. Traditional decentralized solutions tend to ignore the structure of the constraints and lead to degraded performance. We instead develop a decentralized solution that exploits the sparsity structure. Under constant step-size learning, the asymptotic convergence of the proposed algorithm is established in the presence of non-smooth terms, and it occurs at a linear rate in the smooth case. We also examine how the performance of the algorithm is influenced by the sparsity of the constraints. Simulations illustrate the superior performance of the proposed strategy.
\end{abstract}


%\IEEEpeerreviewmaketitle

\section{Introduction}
In many applications such as network utility maximization \cite{palomar2007alternative}, smart grids \cite{halvgaard2016distributed}, basis pursuit \cite{mota2012distributed}, and resource allocation in wireless networks \cite{shen2005adaptive}, a collection of $K$ interconnected agents are coupled through an optimization problem of the following form:
\begin{align}
 \underset{w_1,w_2,\cdots,w_K}{\text{minimize   }}& \quad
  \sum_{k=1}^K J_k(w_k), \ 
 \text{s.t.  } \
 \sum_{k=1}^K B_{k}w_{k}=b, \quad  \label{global-1} 
\end{align}
where $J_k(.)$: $\real^{Q_k} \rightarrow \real$ is a cost function associated with agent $k$ and $w_k \in \real^{Q_k}$ is the variable for the same agent. The matrix $B_{k} \in \real^{S \times Q_{k}}$ is known locally by agent $k$ only, and the vector $b \in \real^{S}$ is known by at least one agent in the network. In this formulation, each agent wants to find its own minimizer, denoted by $w_k^\star$, through interactions with neighboring agents, while satisfying the global coupling constraint. 

In many other applications, the constraint is sparse in the sense that some rows of $B_k$ are zero. For example, in network flow optimization \cite{ahuja1992network}, multitask problems \cite{nassif2017diffusionmulti}, distributed model predictive control \cite{necoara2011parallel}, and optimal power flow \cite{kar2014distributed,dall2013distributed}, the constraint has a special sparse structure. Specifically, each agent $s$ is coupled with its neighboring nodes through an individual affine constraint of the form: 
\eq{
\sum_{k \in \cN_s} B_{s,k} w_k=b_s, \quad \forall \ s=1,\cdots,K \label{constraint-neighbor}
}
where $B_{s,k}\in \real^{S_s \times Q_k}$, $b_{s}\in \real^{S_s}$, and $\cN_s$ denotes the neighborhood of agent $s$ including agent $s$ itself. Note that we can rewrite the constraints \eqref{constraint-neighbor} into a single constraint of the form given in \eqref{global-1} by choosing $B_k$ to be a block column matrix with blocks $\{B_{1,k},\cdots,B_{K,k}\}$ and by setting $B_{s,k}=0$ if $s \notin \cN_k$. However, under {\em decentralized} settings, applying an algorithm that solves \eqref{global-1} directly and ignores the sparsity structure scales badly for large networks and its performance deteriorates as shown in this work. In some other applications (see Example \ref{example:economicdispatch} in Section \ref{section:problem_setup}), unlike \eqref{constraint-neighbor}, the number of constraints is arbitrary, and independent of the number of agents $K$. Moreover, each constraint may include any subset of agents and not only the agents in the neighborhood of some agent. Therefore, a general scalable algorithm that can exploit the sparsity in the constraint set is necessary for large scale networks.
\subsection{Related Works}
 Many distributed/decentralized algorithms have been developed for constraints of the form \eqref{constraint-neighbor}, but for special cases and/or under a different settings from what is considered in this work \cite{giselsson2013accelerated,rostami2017admm,
 kar2014distributed,dall2013distributed,necoara2015linear}. For example, the algorithms developed in \cite{giselsson2013accelerated,rostami2017admm
 ,kar2014distributed,dall2013distributed,necoara2015linear}     require the sharing of {\em primal} variables among neighboring agents and, moreover, the $s-$th constraint is of the form \eqref{constraint-neighbor}, which is limited to agents in the neighborhood of agent $s$. An augmented Lagrangian solution is pursued in \cite{lee2017distributed}, which further requires two hop communications.  All these methods are not directly applicable for the case when the $s$-th constraint involves agents {\em beyond} the neighborhood of agent $s$. Direct extension of these methods to this case would require multi-hop communication, which is costly. Moreover, the settings in these works are different from this work. In these works, the parameters of the $s$-th constraint $\{B_{s,k},b_s\}_{k \in \cN_s}$ are known  by agent $s$. In this work, each agent $s$ is only aware of the constraints matrices multiplying its own vector $w_s$. Moreover, we consider a broader setting with arbitrary number of constraints, and each constraint may involve any subset of agents -- see Section \ref{section:problem_setup}.

The setting in this work is closer to the one considered in \cite{chang2015multi,
chang2014distributed,
chang2016proximal,notarnicola2017constraint,falsone2017dual}. However, these works focused on problems with a single coupling constraint of type \eqref{global-1}, which ignores any sparsity structure. Problem \eqref{global-1} is solved in these references by using dual decomposition methods, which require each agent to maintain a dual variable associated with the constraint. Ignoring any sparsity structure means that each agent will be involved in the entire constraint. By doing so, each agent will maintain a long dual vector to reflect the {\em whole} constraint, and {\em all} agents in the network will have to reach consensus on a longer dual vector.   The work \cite{A1} studied problem \eqref{global-1} for smooth functions with resource constraints (i.e., $\underline{w}_k \leq w_k\leq \overline{w}_k$) and focused on handling the useful case of dynamic and directed graphs.  Note that the matrix $B_k$ in \cite{A1} has a specific structure; but, the solution employed also shares the whole dual variable and neglects any sparsity structure. In other resource allocation problems \cite{ho1980class,xiao2006optimal
,lakshmanan2008decentralized}, all agents are involved in a single constraint of the form \eqref{global-1} with $B_k=I$.

Different from the previously mentioned works, we consider a broader class of coupled affine constraints, where there exist multiple affine constraints and each constraint may involve any connected subset of agents. Our solution requires sharing dual variables only and does not directly share any sensitive primal information, e.g., it does not share the local variables $\{w_k\}$.
Unlike the works \cite{chang2015multi,
chang2014distributed,
chang2016proximal
,notarnicola2017constraint,falsone2017dual,A1}, which solve problem \eqref{global-1} and do not consider the sparsity structure in the constraint, this work exploits the constraint structure. In this way, each agent will only need to maintain the dual variables corresponding to its part of the constraints and not the {\em whole} constraint. Thus, only the agents involved in one particular part will need to agree on the associated dual variables. An algorithm that ignores the sparsity structure scales badly (in terms of communications and memory) as the number of constraints or agents increases.  Moreover, it is theoretically shown in this work that the sparsity in the constraint set influences the performance of the algorithm in terms of convergence rate. Therefore, for large scale networks, it is important to design a scalable algorithm that exploits any sparsity in the constraint.

 In \cite{nassif2017diffusionmulti}, a multi-agent optimization problem is considered with stochastic quadratic costs and an arbitrary number of coupled affine constraints with the assumption that the agents involved in one constraint form a fully connected sub-network. This strong assumption was removed in \cite{hua2017penalty} to handle constraints similar to what is considered in this work albeit with substantially different settings. First, the work \cite{hua2017penalty} considers quadratic costs only, does not handle non-differentiable terms, and their solution solves an approximate penalized problem instead of the original problem. Second, it is assumed that every agent knows all the matrices multiplying the vectors of all other agents involved in the same constraint. For example, for the constraint \eqref{constraint-neighbor}, agent $s$ knows $\{B_{k',k}\}$ for all $k \in  \cN_s$ or $k' \in  \cN_s$. Lastly, the solution method requires every agent to maintain and receive delayed estimates of primal variables $w_k$ from all agents involved in the same constraint through a multi-hop relay protocol. This solution method suffers from high memory and communication burden; thus, it is impractical for large scale networks. 

In network utility maximization problems, a similar formulation appears, albeit with a different distributed framework; it is assumed that the agents (called sources) involved in a constraint are connected through a centralized unit (called link) that handles the constraint coupling these agents  -- see \cite{palomar2007alternative} and references therein. Finally, in \cite{erseghe2012distributed,notarnicola2018distributed} a different ``consensus" formulation is considered where the agents are interested in minimizing  an aggregate cost function where two agents $k$ and $s$ would share similar block vectors  $\{w^k,w^s\}$ if, and only, if they are neighbors, where the notation $w^k$ stands for the block variable shared by the neighbors of agent $k$ so that each $w_k=\mbox{\rm col}\{w^s\}_{s\in {\cal N}_k}$. A more general ``consensus" formulation appears in \cite{alghunaim2019distributed,mota2015localdomains} where the sharing of block entries is not limited to neighboring agents.
\subsection{Main Contributions}
Given the above, we now state the main contributions of this work. A novel low computational decentralized algorithm is developed that exploits the sparsity in the constraints. The developed algorithm handles non-differentiable terms and is shown to converge to the optimal solution for {\em constant} step-sizes.  Furthermore, linear convergence is shown in the absence of non-differentiable terms and an explicit upper bound on the rate of convergence is given. This bound shows the importance of exploiting any constraint sparsity and why not doing so degrades the performance of the designed algorithm.

\noindent {\bf Notation}.  All vectors are column vectors unless otherwise stated. All norms are 2-norms unless otherwise stated. The notation $\|x\|^2_{D}$ denotes the weighted norm $x\tran Dx$ for a positive definite matrix $D$ (or scalar).  The symbol $I_S$ denotes the identity matrix of size $S$ while the symbol $\one_N$ denotes the $N \times 1$ vector with all of its entries equal to one.  We write ${\rm col}\{x_j\}_{j=1}^{N}$ to denote a column vector formed by stacking $x_1, ... , x_N$ on top of each other and $\text{blkdiag}\{X_j\}_{j=1}^{N}$ to denote a block diagonal matrix consisting of diagonal blocks $\{X_j\}$. We let $\text{blkrow}\{X_j\}_{j=1}^{N}=[X_1 \  \cdots \ X_N]$. For the integer set $\cX=\{m_1,m_2,\cdots,m_N\}$, we let $U=[g_{mn}]_{m,n \in \cX}$ denote the $N \times N$ matrix with $(i,j)-$th entry equal to $g_{m_i,m_j}$. The subdifferential $\partial_x f(x)$ of a function $f(.):\real^{M} \rightarrow \real \cup \{\infty\}$ at some $x \in \real^{M}$ is the set of all subgradients:
\eq{
\partial_x f(x) = \{g_x \ | \ g_x\tran(y-x)\leq f(y)-f(x), \forall \ y \in \real^{M}\} \label{subgradients}
}
The proximal operator relative to a function $R(x)$ with step-size $\mu$ is defined by \cite{parikh2014proximal}:
\eq{
{\rm prox}_{\mu R}(x)\define \argmin_u \left(R(u)+{1 \over 2 \mu} \|x-u\|^2 \right)
\label{def_proximal}
} 
\begin{table}[h]
	\centering \caption{A listing of repeatedly used symbols in this work.}
	\begin{tabular}{ | c || l | }
		\hline  \hline
		\cellcolor{gray!25} \bf Symbol &   \multicolumn{1}{c|}{ \bf \cellcolor{gray!25} Description} \\
	 \hline
		$\mathcal{C}_e$ & Sub-network of nodes involved in constraint $e$.  
	 \\ 
		 \hline 
		$N_e$ & The cardinality of the set $\cC_e$.  \\ 
	\hline
	$\cE_k$ & The set of equality constraints indices involving agent $k$.  \\  
			\hline 
		$\sw$ & The vector formed by stacking $\{w_k\}$ over all agents.\\   \hline 
		$\cJ(\sw)$ & The sum of all smooth functions, $\cJ(w)=\sum_{k=1}^K J_k(w_k)$.\\    
	 \hline
		$v^e$ & Dual variable for equality constraint(s) $e$. 
	\\	 \hline  	 		 
		$\{v^e\}_{e \in \cE_k}$ & Collection of all dual variables $v^e$ related to agent $k$. \\
		\hline  	 		 
		$v_k^e$ & Local copy of $v^e$ at agent $k \in \cC_e$.  \\  \hline  	 		 
		$\sy^e$ & Collection of  $v_k^e$ over all $k \in \cC_e$. \\
		\hline  	 		 
		$\sy$ & Collection of $\sy^e$ over all $e$.  \\
				\hline \hline
	\end{tabular}
	\label{table-notation}
\end{table}
\section{Problem Formulation} \label{section:problem_setup}
Consider a network of $K$ agents and assume that the agents are coupled through $E$ affine equality constraint sets. For each constraint set  $e$, we let $\cC_e$ denote the {\em sub-network} of agents involved in this particular constraint(s). We then formulate the following optimization problem:
%
%However, in some cases the structure of the matrices $\{B_k\}$ can be exploited to develop better algorithms. For example, consider the problem:
\begin{align}
 \underset{w_1,\cdots,w_K}{\text{minimize   }}& \quad
   \sum_{k=1}^K J_k(w_k)+R_k(w_k) \label{global-2}  \\
 \text{subject to     }& \quad
 \sum_{k \in \cC_e} \left( B_{e,k}w_{k} -b_{e,k} \right)=0, \quad \forall \ e=1,\cdots,E, \nonumber
\end{align}
where $B_{e,k}\in \real^{S_e \times Q_k}$ and $b_{e,k}\in \real^{S_e}$. The function $J_k(.):\real^{Q_k} \rightarrow \real$ is a smooth function, while $R_k(.): \real^{Q_k} \rightarrow \real \cup \{+\infty\}$ is a convex function possibly non-smooth. For example, $R_k(.)$ could be an indicator function of some local constraints (e.g., $ w_k \geq 0$). These functions are assumed to satisfy the conditions in Assumption \ref{cost-assump} further ahead. It is also assumed that agent $k \in \cC_e$ is only aware of $B_{e,k}$ and $b_{e,k}$. Note that for the special case $E=1$ and $\cC_1=\{1,\cdots,K\}$, problem \eqref{global-2} reduces to \eqref{global-1}.
\begin{assumption} \label{cost-assump}
{\bf (Cost function)}: It is assumed that the aggregate function, $\cJ(\sw)=\sum_{k=1}^K J_k(w_k)$ where $\sw\triangleq{\rm col}\{w_k\}_{k=1}^K$, is a convex differentiable function  with Lipschitz continuous gradient:
\eq{
\|\grad \cJ(\sw)-\grad \cJ(\sw^\bullet)\| \leq \delta \|\sw-\sw^\bullet\|  \label{lipschitz}
}
Moreover, $\cJ(\sw)$ is also strongly convex, namely, it satisfies:
\eq{
(\sw-\sw^\bullet)\tran \grad \cJ(\sw) \geq  \cJ(\sw)- \cJ(\sw^\bullet)+ {\nu \over 2} \|\sw-\sw^\bullet\|^2  \label{stron-convexity}
}
\noindent where $\{\delta,\nu\}$ are strictly positive scalars with  $\delta \geq \nu $. The regularization functions $\{R_k(.)\}$ are assumed to be proper and closed convex functions.
\qd
\end{assumption}
These assumptions are widely employed in the distributed optimization literature and  they are encountered in some practical applications such as distributed model predictive control \cite{necoara2011parallel}, power systems \cite{kar2014distributed}, and data regression problems \cite{chang2015multi}. 
\begin{assumption}
{\bf (Sub-networks)}: \label{assump_connected}
The network of $K$ agents is undirected (i.e., agents can interact in both directions over the edges linking them) and each sub-network ${\cal C}_e$ is connected.
			  \qd
	\end{assumption}
This assumption means that there exists an undirected path between any two agents in each sub-network.  This is automatically satisfied in various applications due to the physical nature of the problem. This is because coupling between agents often occurs for agents that are located close to each other. Applications where this assumption holds include, network flow optimization \cite{ahuja1992network}, optimal power flow \cite{kar2014distributed,dall2013distributed}, and distributed model predictive control \cite{necoara2011parallel} problems. As explained in the introduction, in these problems, the constraints have the form given in equation \eqref{constraint-neighbor}. In this case, each constraint involves only the neighborhood of an agent, so that  ${\cal C}_e={\cal N}_s$ (for $s=e$) and neighborhoods are naturally connected.  Now, more generally, even if some chosen sub-network happens to be disconnected, we can always construct a larger connected sub-network as long as the {\em entire} network is connected -- an explanation of this construction procedure can be found in \cite{alghunaim2017allerton}. The problem of finding this  construction is the well known {\em Steiner tree problem} \cite{hwang1992steiner} and, many decentralized algorithms and heuristics exist to solve it \cite{chalermsook2005simple,bezenvsek2014survey}. We now provide one motivational physical application that also satisfies the two previous assumptions.
%\begin{example}{\rm {\bf (General Exchange Problem)} \label{example:generalexchange}
%This example illustrates the general exchange problem described in \cite{parikh2014proximal}, which fits into \eqref{global-2}. Assume we have $K$ agents and $E$ commodities (or goods). Each commodity $e$ is to be exchanged over a subset of agents participating in the market related to commodity $e$, while minimizing the social welfare. To formalize this, let $\cC_e$ be the subset of agents participating in the exchange of commodity $e$. Moreover, let $w_{e,k} \in \real^{Q'_{e}}$ be the value of commodity $e$ at agent $k \in \cC_e$. If we collect agent $k$  variables $\{w_{k,e} \ | \ e:k \in \cC_e \}$ into an augmented vector $w_k={\rm col}\{w_{e,k}\}_{e:k \in \cC_e} $.
%Then, the general exchange problem is \cite{parikh2014proximal}:
%\begin{align}
% \underset{w_1,\cdots,w_K}{\text{minimize   }}& \quad
%   \sum_{k=1}^K \big(J_k(w_k)+R_k(w_k)\big), \label{global-exchange}  \\
% \text{subject to     }& \quad
% \sum_{k \in \cC_e} w_{e,k}=0, \quad \forall \ e=1,\cdots,E, \nonumber
%\end{align}
%The above formulation occurs in dynamic energy exchange in smart grids applications \cite{kraning2014dynamic} and in the economic dispatch problem \cite{wood2012power} as we illustrate in the next example.
%\qd
%}
%\end{example}
\begin{example}{\rm {\bf (General exchange in smart-grids)} \label{example:economicdispatch}
  For simplicity, we describe the resource management (or economic dispatch) problem in smart grids \cite{wood2012power}  with minimum notation. To begin with, let $P_{G_k}$ and $P_{L_k}$ be the power generation supply and power load demand at node $k$. Moreover, let $P_k={\rm col}\{P_{G_k},P_{L_k}\}$ be a $2 \times 1 $ vector formed by stacking $P_{G_k}$ and $P_{L_k}$.  Then, the resource management problem over a power network consisting of $K$ nodes is \cite{xu2015distributed}: 
 \begin{align}
\scalemath{0.96}{ \underset{\{P_{k}\}}{\text{min}} \
   \sum_{k=1}^K J_k(P_k)+R_k(P_k), \ \text{s.t.} \
 \sum_{k =1}^K \big(P_{G_k}-P_{L_k}\big)=0, } \label{global-smart-grid0}  
\end{align}
where the non-differentiable term $R_k(P_k)$ is the indicator function of some capacity constraints such as positive powers and the maximum power generation. This problem fits into \eqref{global-1} and couples all nodes in a single constraint.  The cost function typically used by power engineers is quadratic and satisfies Assumption \ref{cost-assump} -- see \cite{xu2015distributed,kar2014distributed}. In this formulation, it is assumed that each node is associated with one generator or load with $P_k$ denoting the power generation or demand at that node.  
 Assume now that each node $k$ has multiple generators and/or loads. For example, each generator (or load) can be divided into sub-generators (or sub-loads). Moreover, assume that the power network is divided into $K$ nodes that provide power to $E$ sub-areas. Let $P_{e,G_k}$ and $P_{e,L_k}$ denote the power supply and power load at node $k$ in area $e$ -- see Figure \ref{fig:smart_grid}. In this figure, there are six nodes (agents) and three sub-areas (sub-networks). Each node associates different generators or loads to different sub-areas.
 \begin{figure}[H]
	\centering
	\includegraphics[width=0.5\textwidth]{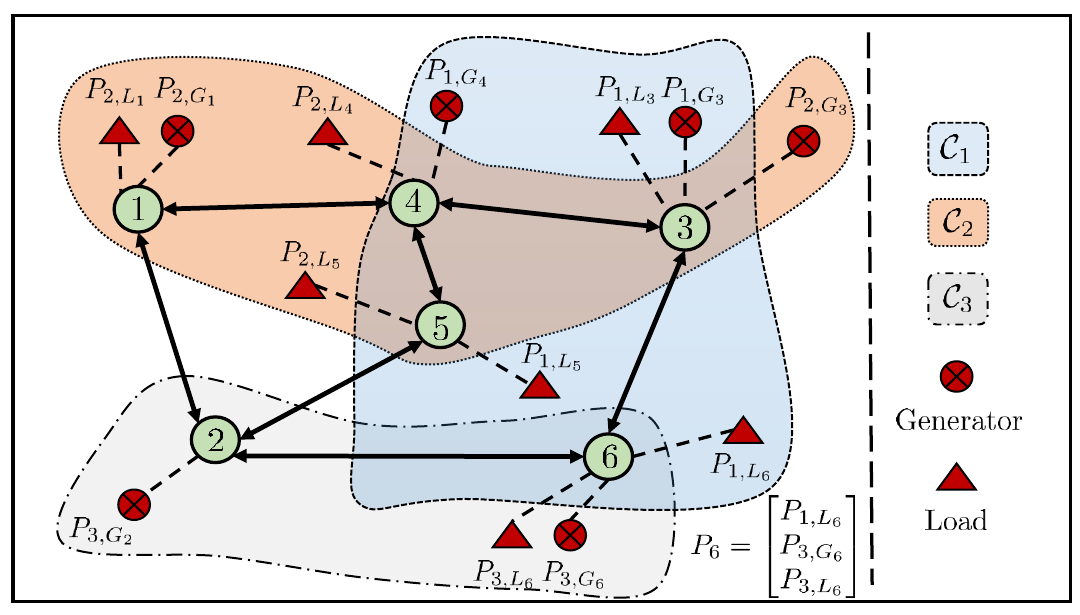}
	\caption{\small An illustration for Example \ref{example:economicdispatch}. In this illustration, there are $E=3$ areas and $K=6$ agents.}
	\label{fig:smart_grid}
\end{figure} 
\noindent If we let $\cC_e$ denote the nodes that are involved in area $e$ and $P_k$ to be the augmented vector $P_k={\rm col}\{P_{e,G_k},P_{e,L_k}\}_{e:k \in \cC_e}$, which collects all local variables $\{P_{e,G_k},P_{e,L_k}\}$ over all areas that agent $k$ belongs to. Then, we formulate the following more general problem:
\begin{align}
 \underset{\{P_{k}\}}{\text{minimize   }}& \quad
   \sum_{k=1}^K J_k(P_k)+R_k(P_k) \label{global-smart-grid}  \\
 \text{subject to     }& \quad
 \sum_{k \in \cC_e}   (P_{e,G_k}-P_{e,L_k})=0, \quad \forall \ e=1,\cdots,E \nonumber
\end{align}
   This formulation fits into the problem of  dynamic energy exchange in smart grids applications \cite{kraning2014dynamic} where each area satisfies Assumption \ref{assump_connected}. It can also be motivated as follows. Assume each sub-area represents some city. Then, problem \eqref{global-smart-grid} is useful when the transmission losses are costly in some parts of an area, which may require power generation from neighboring power networks. It is also useful when there are maintenance to some generators or lines causing high demands in some areas, which requires the need of extra generators from adjacent power networks.
\qd
}
\end{example}
\begin{figure*}[t]
	\centering
	\includegraphics[scale=0.5]{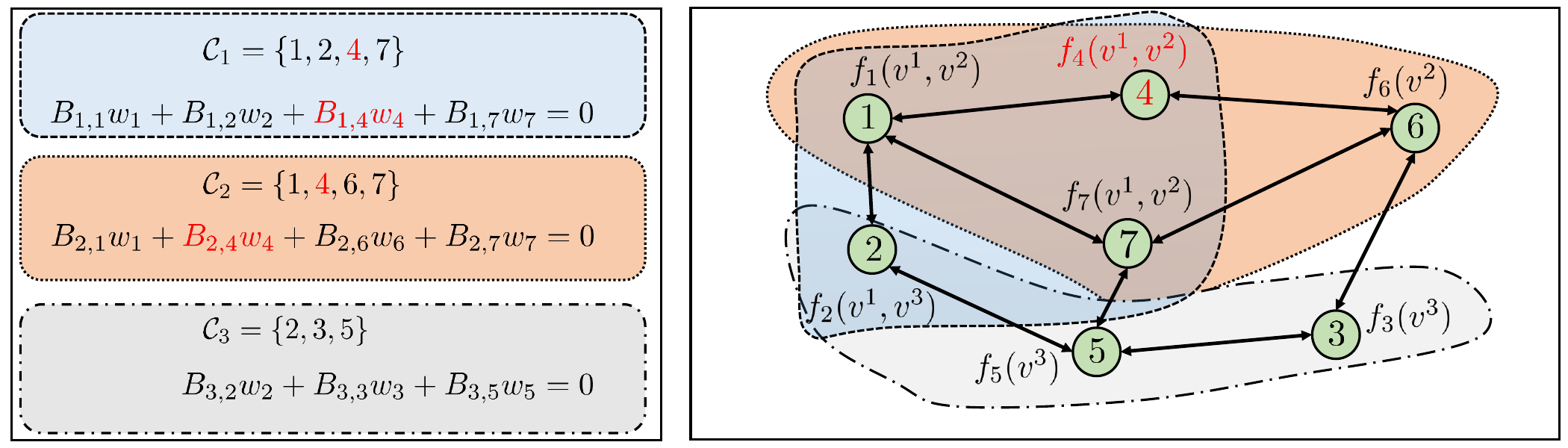}
	\caption{\small An example of  problem \eqref{dual_problem} with three sub-networks with  agent $4$ involved in the equality constraints for sub-networks ${\cal C}_1$ and ${\cal C}_2$.}
	\label{fig:network-dual2}
\end{figure*}
\section{Algorithm Development}
 In this section, we will derive  our algorithm and introduce some important symbols,  which are necessary for algorithm description and later analysis. To do so, we start by introducing the Lagrangian function of \eqref{global-2}:
\eq{
\hspace{-2mm}  \scalemath{0.95}{ \cL\big(\sw,\{v^e\}_{e=1}^E\big) \hspace{-0.75mm}=\hspace{-0.75mm} \sum_{k=1}^K J'_k(w_k) \hspace{-0.5mm}+ \hspace{-0.5mm} \sum_{e=1}^E (v^e)\tran \hspace{-0.5mm} \bigg(\hspace{-0.75mm} \sum_{k \in \cC_e} \hspace{-0.75mm} B_{e,k}w_{k} \hspace{-0.5mm} -\hspace{-0.5mm} b_{e,k} \hspace{-1mm} \bigg) } \hspace{-2mm}  \label{Lagrangian}
}
where  $J'_k(w_k) \define J_k(w_k)+R_k(w_k)$, and $v^e \in \real^{S_e}$ denotes the dual variable associated with the $e$-th constraint.  To facilitate the development of the algorithm we rewrite \eqref{Lagrangian} as a sum of local Lagrangian terms. To do so, we need to introduce the set $\cE_k$, which denotes the set of equality constraints that agent $k$ is involved in (e.g., if agent $k$ is involved in equality constraints one and three, then $\cE_{k}=\{1,3\}$).  From the definition of $\cE_k$ and $\cC_e$, we have 
\eq{
\cC_e=\{k \ | \ e\in \cE_k\}, \quad \cE_k=\{e \ | \ k\in \cC_e\} \label{equality_cal_E_C} 
}
Using this notation, the second term on the right hand side of \eqref{Lagrangian} can be rewritten as a sum over all agents as follows: let $\overline{B}_{e,k}=B_{e,k}$ if $k \in \cC_e$ (or $e \in \cE_k$) and zero otherwise and, likewise, for $\overline{b}_{e,k}$. then it holds that:
\eq{
\scalemath{0.95}{\sum_{e=1}^E  \sum_{k \in \cC_e} (v^e)\tran \left(B_{e,k}w_{k}-b_{e,k}\right) } 
 &=\scalemath{0.95}{\sum_{e=1}^E  \sum_{k=1}^K (v^e)\tran \left(\overline{B}_{e,k}w_{k}-\overline{b}_{e,k}\right)} \nonumber \\
 &  =
 \scalemath{0.95}{\sum_{k =1}^K \sum_{e \in \cE_k} (v^e)\tran \left(B_{e,k}w_{k}-b_{e,k}\right)} 
\nonumber}
where in the last step we switched the order of summation and used the fact that $k \in \cC_e$ if, and only, if $e \in \cE_k$.
 Therefore, if we let $\{v^e\}_{e \in \cE_k}$ denote the collection of dual variables related to agent $k$,
then using the previous equation we can rewrite \eqref{Lagrangian} as a sum of local terms as follows: 
\eq{
\cL\big(\sw,\{v^e\}_{e=1}^E\big) = \sum_{k=1}^K L_k\big(w_k,\{v^e\}_{e \in \cE_k}\big) \label{Lagrangian2}}
where
\eq{
\scalemath{0.95}{ L_k\big(w_k,\{v^e\}_{e \in \cE_k}\big) \hspace{-1mm} \define \hspace{-1mm} J'_k(w_k) \hspace{-0.5mm}+\hspace{-0.5mm} \sum_{e \in \cE_k} (v^e)\tran \big(  B_{e,k}w_{k}- b_{e,k} \big) }
\label{local_term}
}
 is the local term for agent $k$.  We are therefore interested in finding the minimizer of \eqref{global-2} through the equivalent solution of the saddle point problem:
\eq{
\min_{\ssw} \max_{\{v^e\}_{e=1}^E} \quad \cL(\sw,\{v^e\}_{e=1}^E) \label{max_min}
}
\begin{assumption}{\bf (Strong duality)} \label{assum_strong_dual}
A solution exists for problem \eqref{max_min} and strong duality holds. \qd
\end{assumption}
 Since our problem \eqref{global-2} is convex with affine constraints only, then Slater's condition is satisfied and strong duality holds \cite[Section~5.2.3]{boyd2004convex}, which ensures that the solution of \eqref{max_min} coincides with the solution of \eqref{global-2}. We denote an optimal solution pair of \eqref{max_min} by $\sw^\star={\rm col}\{w_k^\star\}_{k=1}^K$ and $\{v^{e,\star}\}$. From Assumption \eqref{cost-assump}, $\sw^\star$ is unique, but $\{v^{e,\star}\}$ are not necessarily unique.  To derive our algorithm, which solves the saddle point problem \eqref{max_min}, we will now relate the dual problem to the one considered in our previous work \cite{alghunaim2017allerton} and explain how the dual variables are partially shared across the agents, which is important for our derivation. 
\subsection{Dual Problem} 
 Note that the Lagrangian \eqref{Lagrangian2} is separable in the variables $\{w_k\}$. Thus, the dual problem is (we are reversing the $\min$ and $\max$ operations by negating the function) \cite{boyd2004convex}:
  \begin{align}
 \underset{v^1,\cdots,v^E}{\text{minimize   }}& \quad
  -\sum_{k=1}^K f_k\big(\{v^e\}_{e \in \cE_k} \big) \label{dual_problem}  
\end{align}
where\footnote{Technically inf should be used instead of min in \eqref{Lmin}, however, to avoid confusion we use min.}
  \eq{
 f_k\big(\{v^e\}_{e \in \cE_k} \big)\define \min_{w_{k}} L_k\big(w_k,\{v^e\}_{e \in \cE_k} \big)
 \label{Lmin} }
 Figure \ref{fig:network-dual2} illustrates how the dual variables $\{v^e\}_{e=1}^E$ are shared across agents participating in the same constraint. For example, agent $k=4$ in Figure \ref{fig:network-dual2} is part of two sub-networks, ${\cal C}_1$ and ${\cal C}_2$; it is therefore part of two equality constraints and will be influenced by their respective dual variables, denoted by $v^1$ and $v^2$. Similarly, for the other agents in the network. Problem \eqref{dual_problem} is of the form considered in \cite{alghunaim2017allerton}:  it involves minimizing the aggregate sum of cost functions $f_k\big(\{v^e\}_{e \in \cE_k} \big)$ where the arguments $\{v^e\}_{e \in \cE_k}$ among different agents can share block entries as illustrated in Fig. \ref{fig:network-dual2} . The main difference here, however, is that the costs $f_k\big(\{v^e\}_{e \in \cE_k} \big)$ do not admit a closed form expression in general and are instead defined by \eqref{Lmin}, i.e., in this work we are actually dealing with the more challenging decentralized {\em saddle point} problem and not with a decentralized minimization problem as was the case in \cite{alghunaim2017allerton}. Thus, more is needed to arrive at the solution of \eqref{max_min}, as we explain later. 
\subsection{Combination coefficients}
  To proceed from here and for the algorithm description, we introduce combination coefficients for the edges in ${\cal C}_{e}$ denoted by $\{a_{e,sk}\}_{s,k \in \cC_e}$;   $a_{e,sk}$ refers to the coefficient used to scale data moving from agent $s$ to agent $k$ in subnetwork ${\cal C}_e$ with $a_{e,sk}=0$ if $s \notin \cN_k \cap \cC_e$. We collect these coefficients into the combination matrix
 \eq{
 A_e\define [a_{e,sk}]_{s,k \in \cC_e} \in \real^{N_e \times N_e} \label{A_e_combination}
 } 
where $N_e$ denotes the number of agents involved in equality $e$. The matrix $A_e$ is assumed to be symmetric and doubly-stochastic. We also require $A_e$ to be primitive, meaning that  there exists an integer $j$ such that the entries of the matrix $A_e^j$ are all positive.  One way to meet these conditions is to choose weights satisfying
\begin{subnumcases}{\label{primitive_assump}}
\sum_{s \in \cC_e} a_{e,sk}=1, \quad \quad \sum_{k \in \cC_e} a_{e,sk}=1 \label{com_coeff_1} \\
a_{e,sk}>0 \ \ {\rm for} \ \ s \in \cN_k \cap \cC_e
\label{com_coeff_primitve} 
\end{subnumcases}   
with $a_{e,sk}=0$ if $s \notin \cN_k \cap \cC_e $.   Under Assumption \eqref{assump_connected} many rules exists to choose such weights in a decentralized way -- see \cite[Ch. 14]{sayed2014nowbook}. We are now ready to derive our algorithm. 
\subsection{Dual Coupled Diffusion}
Using the combination matrix $A_e$, it was shown in \cite{alghunaim2017allerton} that  problem \eqref{dual_problem} can be solved by using the following coupled diffusion algorithm. Set $v^e_{k,-1}=\psi^e_{k,-1}$ to arbitrary values. For each $k$ and $e\in {\cal E}_k$ repeat for  $i\geq 1$:
\begin{subequations}
\eq{
\psi^e_{k,i}&=v^e_{k,i-1}+\mu_v \grad_{v^e} f_k\big(\{v^e_{k,i-1}\}_{e \in \cE_k}\big) \label{exact-diff1}\\
\phi^e_{k,i}&=\psi^e_{k,i}+v^e_{k,i-1}-\psi^e_{k,i-1} \label{exact-diff2}\\
v^e_{k,i}&=\sum_{s \in \cN_k \cap \cC_e} \bar{a}_{e,sk}\phi^e_{s,i}
\label{exact-diff3}
} 
\end{subequations}
where $v^e_{k,i}$ is the estimate for $v^e$ at agent $k$, $\mu_v>0$ is a step-size parameter, and  $\{\psi^e_{k,i},\phi^e_{k,i}\}$ are auxiliary vectors used to find $v^e_{k,i}$.  The coefficients $\{\bar{a}_{e,sk}\}$ are the entries of the matrix $\bar{A_e}$ defined as follows:
 \eq{
  \bar{A}_e = [\bar{a}_{e,sk}]_{s,k \in \cC_e} \define 0.5(I_{N_e}+A_e) \label{barA_e_combination}
 } 
 \begin{remark}[\sc Combination Weights]{\em 
Note that since $A_e$ is primitive, symmetric, and doubly stochastic, it holds that the eigenvalues of the matrix $A_e$ are in $(-1,1]$ -- see \cite[Lemma F.4]{sayed2014nowbook}. Equation \eqref{barA_e_combination} implies that the eigenvalues of the matrix  $\bar{A}_e$ are in $(0,1]$.   The recent preprint \cite{li2019linear}  studies \eqref{exact-diff1}--\eqref{exact-diff3} for one variable $v^1$ ($E=1$ and $\cC_1=\{1,\cdots,K\}$). It is shown in \cite{li2019linear}  that the eigenvalues of  $A_1$ can be relaxed to be in $(-{5 \over 3},1]$ so that the eigenvalues of $\bar{A}_1$ are in $(-{1 \over 3},1]$.     }\qd
\end{remark} 
If the functions $\{ f_k(.)\}$ are known and are differentiable, then each agent could run \eqref{exact-diff1}--\eqref{exact-diff3} to converge to its corresponding optimal dual variable, which in turn could be used to find the local minimizer $w_k^\star$ by solving $\min_{w_k} L_k\big(w_k,\{v^{e,\star}\}_{e \in \cE_k}\big)$.  However, this approach is not always possible because the local dual function $ f_k\big(\{v^e\}_{e \in \cE_k} \big)$ does not generally admit a closed form expression. Moreover, this method involves two time scales: one for finding the dual and the other for finding the primal. Therefore, to solve \eqref{max_min} we propose to employ a decentralized version of the centralized dual-ascent construction \cite{boyd2011admm} combined with a proximal gradient descent step. Specifically, recall first that the  dual-ascent method updates the primal variable $w_k$ at each iteration $i$ as follows:
\eq{
w_{k,i}&= \argmin_{w_k} L_k(w_k,\{v^e_{i-1}\}_{e \in \cE_k}), \ \forall \ k \label{centralized_primal}
%y^e_i&=y^e_{i-1}+\mu_v \grad_{y^e} L (w_i,y_{k,i-1}), \forall \ e \label{centralized_dual}
}
 Note that this minimization step, which need to be solved at each iteration, can be costly in terms of computation unless a closed form solution exists, which is not the case in general. Therefore, we approximate \eqref{centralized_primal} by a proximal gradient descent step to arrive at what we shall refer to as the {\em dual coupled diffusion} algorithm \eqref{final-algorithm_jk=1}. At each time instant $i$, each agent $k$ first performs a proximal gradient descent step \eqref{proximal-grad} for the primal variable with step-size $\mu_w>0$. Then, for each dual-ascent step, the coupled diffusion \eqref{grad_ascent_diff}--\eqref{grad_comb_diff} are applied where step \eqref{grad_ascent_diff} is obtained by using  $\grad_{v^e} L_k(w_{k,i},\{v^e_{k}\}_{e \in \cE_k})$ to approximate the gradient at the minimum value in \eqref{Lmin}. Note that only step \eqref{grad_comb_diff} requires sharing dual variables with the neighbors that are involved in similar constraints. We remark that
Algorithm \eqref{final-algorithm_jk=1} can be potentially used for directed network if combined with the push-sum technique from \cite{kempe2003gossip} such that  the dual iterates are corrected by dividing them by scalar as in \cite{kempe2003gossip}.  The push-sum technique have been utilized before for distributed optimization algorithms -- see for example \cite{nedic2017achieving}.

 To analyze algorithm  \eqref{final-algorithm_jk=1} and show that it converges to an optimal solution of \eqref{max_min},  we will rewrite it in a compact network form, which facilitates its analysis.
\begin{algorithm}[t] 
\caption*{\textrm{\bf{Algorithm}} (\sc Dual Coupled Diffusion)}
{\bf Setting}:  Choose step-sizes $\mu_w>0$ and $\mu_v>0$. Let $v^e_{k,-1}=\psi^e_{k,-1}$ and $w_{k,-1}$ arbitrary. 

{\bf For every agent $k$, repeat for $i \geq 0$:}
\begin{subequations} \label{final-algorithm_jk=1}
\eq{
w_{k,i}&=\scalemath{0.95}{\underset{\mu_{w} R_{k}}{{\rm prox}} \big(w_{k,i-1}  -\mu_w  \grad J_k(w_{k,i-1})
 - \mu_w \sum_{e \in \cE_k} B_{e,k}\tran  v^e_{k,i-1} \big)} \label{proximal-grad}}
\hspace{10mm} { For all $e \in \cE_k$}:
 \eq{
\psi^e_{k,i}&=v^e_{k,i-1}+\mu_v  \big( B_{e,k}w_{k,i}- b_{e,k} \big) \label{grad_ascent_diff} \\
\phi^e_{k,i}&=\psi^e_{k,i}+v^e_{k,i-1}-\psi^e_{k,i-1} \label{grad_correct_diff} \\
v^e_{k,i}&=\sum_{s \in \cN_k \cap \cC_e} \bar{a}_{e,sk}\phi^e_{s,i} 
\label{grad_comb_diff}}
\end{subequations}
\end{algorithm}
\section{Network Recursion}
\begin{figure}[H]
	\centering
	\includegraphics[scale=0.275]{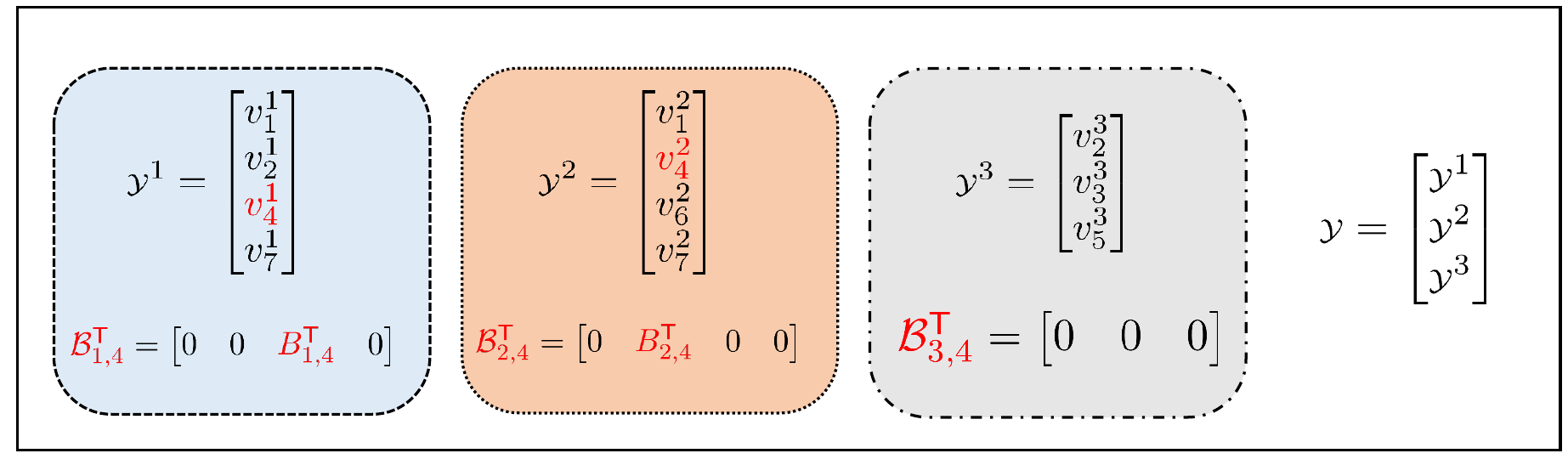}
	\caption{\footnotesize An illustration of constructions \eqref{cal_Y_defin_cluster} and \eqref{cal_Y_defin} for the network in Figure \ref{fig:network-dual2} as well as construction \eqref{calBblkrow} for agent $k=4$ in that network.}
	\label{fig:constructionB}
\end{figure}

  We start by stacking the dual estimates within each cluster and then stacking over all the clusters. This will allow us to rewrite the dual steps \eqref{grad_ascent_diff}--\eqref{grad_comb_diff} in a form that enables us to see the affect of each sub-network in our analysis. Thus, we introduce the sub-network vector that collects the dual estimates $v_{k,i}^e$ over the agents in $\cC_e$:
\eq{
\sy^e_{i}& \define {\rm col}\{v^e_{k,i}\}_{k \in \cC_e}  \in \real^{N_e S_e}, \label{cal_Y_defin_cluster}
}
and the global network vector that collects $\sy^e_{i}$ over all $e$: 
\eq{
\sy_i &\define {\rm col}\{\sy^e_{i}\}_{e=1}^E \label{cal_Y_defin}  
}
We also repeat a similar construction for the quantities:
\eq{
b_e& \define {\rm col}\{b_{e,k}\}_{k \in \cC_e}, \quad &b&\define {\rm col}\left\{b_e \right\}_{e=1}^E \label{bedef} \\
 \bar{\cA}_e &\define \bar{A}_e \otimes I_{S_e}, \quad
 &\bar{\cA} &\define {\rm blkdiag}\{\bar{\cA}_e\}_{e=1}^E \label{cal_bar_A_defin} 
}
where $\bar{A}_e={1\over2}(I_{N_e}+A_e)$ introduced in \eqref{barA_e_combination}. For the networked representation of the primal update \eqref{proximal-grad}, we introduce the network quantities:
\eq{
\sw_{i} &\define {\rm col}\{w_{k,i}\}_{k=1}^K \\
\cR(\sw) &\define \sum_{k=1}^K R_k(w_k) \\ \grad\cJ(\sw_{i}) &\define {\rm col}\{\grad J_k(w_{k,i})\}_{k=1}^K \label{gradcalJ} }
 We also need to represent the term  $\sum_{e \in \cE_k} B_{e,k}\tran  v^e_{k,i-1}$ in terms of the network quantity $\sy_{i-1}$ defined in \eqref{cal_Y_defin}.  To do that we first rewrite each term $B_{e,k}\tran  v^e_{k,i-1}$ in terms of the sub-network vector $\sy_{i-1}^e$. This can be simply done by introducing the $1 \times N_e$ block row matrix $\cB_{ek}\tran$ of similar block structure as $\sy_{i-1}^e$ such that $\cB_{ek}\tran  \sy^e_{i-1}=B_{e,k}\tran  y^e_{k,i-1}$ if $k \in \cC_e$ and zero otherwise -- Figure \ref{fig:constructionB} illustrates this construction. This construction can be represented by:
  \begin{subequations} \label{calBblkrow}
  \eq{
  \cB\tran_{ek} &= {\rm blkrow}\{B\tran_{e,kk'}\}_{k' \in \cC_{e}} \\
 B\tran_{e,kk'}&\define \begin{cases}
\begin{aligned}
   &B_{e,k}\tran
 ,&\  &\text{if} \ \ k \in \cC_e \ ,k=k' \\
&0_{Q_k,S_e},&\ \ &\text{otherwise}
\end{aligned}
\end{cases}
  }
  \end{subequations}
   Thus, we have $
\sum_{e \in \cE_k} B_{e,k}\tran  v^e_{k,i-1} =\sum_{e =1}^E \cB_{ek}\tran  \sy^e_{i-1}$.
\noindent  If we let
  \eq{
\cB &\define \begin{bmatrix}
\cB_{11} & \cdots & \cB_{1K} \\
\vdots & & \vdots \\
\cB_{E1} & \cdots & \cB_{EK} 
\end{bmatrix} \label{calBtran}
}
then algorithm  \eqref{final-algorithm_jk=1} can be rewritten compactly as follows:
\begin{subequations}
\eq{
\sw_{i}&=\scalemath{1}{\underset{\mu_{w} \cR}{{\rm prox}} \big(\sw_{i-1}-\mu_w \grad\cJ(\sw_{i-1}) -\mu_w \cB\tran \sy_{i-1}\big) }  \label{primal-network_w} \\
\sy_{i}&=\bar{\cA} \bigg(2\sy_{i-1}-\sy_{i-2}+\mu_v   \cB (\sw_{i}- \sw_{i-1}) \bigg)  \label{networkrecursion2}
}
\end{subequations}
for $i \geq 1$ with initialization:
\eq{
\sy_{0}&= \sy_{-1}+\mu_v  ( \cB \sw_{0}-  b)  \label{initial_y}
}
 Notice that step \eqref{networkrecursion2} depends on the two previous estimates; thus it is tedious to analyze directly. Therefore, to facilitate our analysis we will rewrite it in an equivalent form. To do that, we let:
\eq{
\cA = {\rm blkdiag}\{A_e \otimes I_{S_e}\}_{e=1}^E \label{calA}
}
and introduce the singular value (or eigenvalue for symmetric matrices) decomposition \cite{laub2004}:
\eq{
0.5(I_{N}-\cA )= \begin{bmatrix}
\cU_1 &\cU_2
\end{bmatrix} \begin{bmatrix}
\Sigma & 0 \\
0 & 0
\end{bmatrix} \begin{bmatrix}
\cU_1\tran \\ \cU_2\tran
\end{bmatrix} = \cU_1 \Sigma \cU_1\tran
\label{I-Adecomposition}}
where $N=\sum_{e=1}^E N_e S_e$, $\cU_1 \in \real^{N\times r}$, $\cU_2 \in \real^{N\times (N-r)}$, and $
\Sigma={\rm diag}\{\lambda_{j}\}_{j=1}^r$
 with $ \lambda_r \leq \cdots \leq \lambda_{1}$ denoting the non-zero eigenvalues of the matrix $0.5(I-\cA)$. Using an approach similar to the one used in \cite{li2017nids}, we can rewrite \eqref{networkrecursion2} equivalently as follows --- see Appendix \ref{Appendix:dual-represe}:
\begin{subequations}
\eq{
\hspace{-1mm}\scalemath{0.95}{ \ssx_{i}}&= \scalemath{0.95}{\ssx_{i-1} - {1 \over \mu_v} \cU_1\tran \big(\sy_{i-1}+\mu_v ( \cB \sw_{i}-  b)  +\mu_v \cU_1\Sigma \ssx_{i-1} \big) }  \label{dual1-network} \hspace{-1mm} \\
\hspace{-1mm} \scalemath{0.95}{ \sy_{i}}&=\scalemath{0.95}{\sy_{i-1}+\mu_v \big( \cB \sw_{i}-  b\big)+\mu_v \cU_1 \Sigma \ssx_i} \hspace{-1mm} \label{dual2-network}
}
\end{subequations}
for $i\geq 1$, where we introduced a new sequence $\ssx_{i}$
 with $\ssx_{0}=0$.  Note that since $A_e$ is primitive, symmetric, and doubly stochastic, it holds that the eigenvalues of the matrix $A_e$ are in $(-1,1]$ -- see \cite[Lemma F.4]{sayed2014nowbook}.
 Thus, from the block structure of $\cA$ in \eqref{calA}, the eigenvalues of the matrix $0.5(I-\cA)$ are in $[0,1)$. Therefore, the non-zero eigenvalues are positive and satisfy:
 \eq{
 0 < \lambda_r  \leq \cdots \leq  \lambda_{1} < 1\label{eigenvalues_sigma}
 }
 This property is useful for our convergence analysis.
 
\section{Convergence Results}
In this section, we give the Lemmas leading to the main convergence results.
The following auxiliary result is proven in \cite{yuan2019exactdiffI}.
\begin{lemma}\label{null_I-A}
	For any $N \times N$  symmetric and doubly stochastic matrix $A$, it holds that $I_N-A$ is symmetric and positive semi-definite. If in addition $A$ is primitive and we let $\sa=A \otimes I_{M}$,
	then, for any block vector $\ssz=\text{\em col}\{z^1,...,z^{N}\}$ in the nullspace of $I -\sa$ with entries $z^n \in \mathbb{R}^{M}$ it holds that:
	\eq{
		(I -\sa)\ssz=0  \iff z^1=z^2=...=z^{N}
		\label{nullspace}
	} \qd
\end{lemma}	
\noindent  Lemma \ref{null_I-A} will be used in the proof of the next Lemma to show that consensus is reached at the optimality conditions. 
\begin{lemma} {\bf (Optimality condition)} \label{lemma_optimiality}
If there exists a point $(\sw^\star,\sy^\star,\ssx^\star)$ and a subgradient $g^\star \in \partial_{\ssw} \cR(\sw^\star)$  such that:
\begin{subequations} \label{optimality_conditions}
\eq{
\grad\cJ(\sw^\star) +g^\star + \cB\tran \sy^\star&=0 \label{optimality1} \\
\cU_1\tran \sy^\star&=0 \label{optimality3} \\
 (\cB\sw^\star-b)+ \cU_1 \Sigma \ssx^\star&=0 \label{optimality2} 
}
\end{subequations}
Then, it holds that $
v^{e,\star}_{k}=v^{e,\star}$ $\forall \ k \in \mathcal{C}_{e}$
where $(\sw^\star,v^{1,\star},\cdots,v^{e,\star})$ is a saddle point for the Lagrangian \eqref{Lagrangian}.
\end{lemma}
{\bf Proof:}
   A similar argument appears in the conference version \cite[Lemma~2]{alghunaim2018cdc} except for the addition of sub-gradient terms into the argument. Using the block structure of $\grad \cJ(.)$ and $\cB$ in \eqref{gradcalJ} and \eqref{calBblkrow}--\eqref{calBtran}, we can  expand \eqref{optimality1} into its components to get:
\eq{
\grad_{w_k} J_{k}(w_k^\star)+g_k^\star + \sum_{e \in \cE_k} B_{e,k}\tran v_k^{e,\star}=0, \quad \forall \ k \label{KKT1}
}
where $g_k^\star  \in \partial_{w_k} R_k(w_k^\star)$.  From the fact $\cU_1\tran \cU_1=I$ and $\Sigma >0$, condition \eqref{optimality3} is equivalent to:
\eq{
\scalemath{0.96}{\cU_1\tran \sy^\star=0 \iff \cU_1 \Sigma \cU_1\tran \sy^\star=0 \iff {1\over 2}(I-\cA)  \sy^\star=0}
}
Therefore, from \eqref{nullspace}, and the block structure of $\cA$ in \eqref{calA}, condition \eqref{optimality3} gives:
\eq{
v^{e,\star}_{k}=v^{e,\star}_s=v^{e,\star}, \quad \forall \ k,s \in \cC_e \label{consensus_proof_lemma}
}
for some $v^{e,\star}$.  Hence, condition \eqref{KKT1} satisfies the first optimality condition for problem \eqref{global-2} -- see \cite{boyd2004convex}.  Now, let $\cZ={\rm blkdiag}\{\one_{N_{e}} \otimes I_{S_{e}}\}_{e=1}^E$. Multiplying equation \eqref{optimality2} on the left by $\cZ\tran$ gives:
\eq{
0 = \cZ\tran (\cB\sw^\star-b)+ \underbrace{\cZ\tran \cU_1}_{=0}  \Sigma \ssx^\star \overset{(a)}{=} \cZ\tran (\cB\sw^\star-b) \label{kkk}
}
where step (a) holds because because from \eqref{nullspace}, $\cZ$ is in the nullspace of $I-\cA$ and thus also in the nullspace of $\cU_1\tran$ \cite[Equation (51)]{alghunaim2018cdc}. Using the block structure of $\cB$ and $b$ in \eqref{calBblkrow}--\eqref{calBtran} and \eqref{bedef}, we can also expand \eqref{kkk} into its components to get:
\eq{
 (\one_{N_{e}}\tran \otimes I_{S_{e}}) \bigg(\sum_{k=1}^K  (\cB_{ek} w_k^\star)  -  b_e \bigg)  &=\sum_{k \in \cC_e}  \left( B_{e,k} w_k^\star  -b_{e,k} \right) \nnb
 & =0
\label{KKT2}
}
for all $e$ since
\eq{
\hspace{-2mm} \scalemath{0.975}{ \cB_{ek}\tran (\one_{N_{e}} \otimes I_{S_{e}}) }
 &= \scalemath{0.975}{ \sum_{k' \in \cC_e} B_{e,kk'}\tran = \begin{cases}
\begin{aligned}
 &  B_{e,k}\tran
 ,\  &\text{if } k \in \cC_e  \\
&0, &\text{otherwise}
\end{aligned}
\end{cases}
}
}
 Equation \eqref{KKT2} is the second optimality condition for problem \eqref{global-2} and, thus, $(\sw^\star,v^{1,\star},\cdots,v^{e,\star})$ is an optimal point for \eqref{max_min} \cite{boyd2004convex}. 
 \qd
\begin{remark}[\sc Existence and uniqueness]{\em \label{remark_uniqueness}
 Note that there exists a point $(\sw^\star,\sy^\star,\ssx^\star)$ that satisfies the optimality conditions  \eqref{optimality_conditions}. Specifically, if $\sw^\star={\rm col}\{w_k^\star\}$ and $\sy^\star={\rm col}\{\one_{N_e} \otimes v^{e,\star}\}_{e=1}^E$, where $(\sw^\star,v^{1,\star},\cdots,v^{E,\star})$ is an optimal solution of the saddle point problem \eqref{max_min},  then, it can be easily verified that conditions \eqref{optimality1}--\eqref{optimality3} are satisfied. Now, by following an argument similar to the one used in \cite[Lemma ~3]{yuan2019exactdiffII}, it can be shown that there exists an $\ssx^\star$ such that \eqref{optimality2} holds; moreover, there exists a unique $\ssx^\star$ in the range space of $\cU_1\tran$. Now, we know from strong convexity that $\sw^\star$ is unique. Thus, from \eqref{optimality1}, the dual point $\sy^\star$ is unique if the matrix $\cB$ has full row rank. Under this condition and in the absence of non-smooth terms, we will show that our algorithm converges linearly to this unique point -- see Theorem \ref{coro-linear-conv} .}\qd
\end{remark}
\begin{remark}{\em 
The analysis technique used in this work is not related to the techniques used in  \cite{yuan2019exactdiffII,alghunaim2017allerton}. Note that this work deals with a {\em non-smooth saddle-point} problem where the dual variables are shared across agents, while the works \cite{yuan2019exactdiffII,alghunaim2017allerton} deal with {\em smooth minimization} problems with a shared primal variable and twice-differentiable functions. }\qd
\end{remark}
  We will now show that the equivalent network recursions \eqref{primal-network_w} and \eqref{dual1-network}--\eqref{dual2-network} of the proposed algorithm converge to a point that satisfies the optimality conditions given in Lemma \ref{lemma_optimiality}.  To give the convergence results, we introduce the error vectors:
\eq{
\tsw_i \define \sw^\star-\sw_i, \quad \tsx_i \define \ssx^\star-\ssx_{i} \quad \tsy_i &\define \sy^\star-\sy_{i}
}
and the diagonal matrix:
\eq{
\cD \define \mu_v(\Sigma-\Sigma^2) >0 \label{cal_D}
}
where $\Sigma$ was introduced in \eqref{I-Adecomposition}. Note that $\cD$ is positive definite because of \eqref{eigenvalues_sigma}.
\begin{lemma}{\bf (Primal-dual bound)}:
Suppose Assumptions  \ref{cost-assump}-\ref{assum_strong_dual} hold, then:
\eq{
\|\tsw_i\|^2-\|\tsw_{i-1}\|^2 &\leq -\big(1-\mu_w (2\delta-\nu)\big)\|\sw_i-\sw_{i-1}\|^2 \nonumber \\
& \quad - \mu_w \nu \left(\|\tsw_{i-1}\|^2+\|\tsw_i\|^2 \right)    \nonumber \\
& \quad -2\mu_w (\sy_{i-1}-\sy^\star)\tran \cB(\sw_{i}-\sw^\star) \label{werrrrr}
}
and
\eq{
&\|\tsy_{i}\|^2_{\mu_{v}^{-1}}+\|\tsx_{i}\|_\cD^2-\|\tsy_{i-1}\|_{\mu_{v}^{-1}}^2-\|\tsx_{i-1}\|_\cD^2   \nonumber \\
&=  -\|\ssx_{i}-\ssx_{i-1}\|_\cD^2-\|\Sigma \tsx_i \|_{\mu_{v}}^2+ \|\cB \tsw_i\|_{\mu_{v}}^2  \nonumber \\
&\quad + 2(\sy_{i-1}-\sy^\star)\tran \cB (\sw_{i}- \sw^\star)  
 \label{error3a}
}
where $(\sw^\star,\sy^\star,\ssx^\star)$ satisfy the optimality conditions given in Lemma \ref{lemma_optimiality}.
\end{lemma}
{\bf Proof:}  See Appendices \ref{Appendix:Primal_error} and \ref{Appendix:dual_error}. \qd

\noindent The previous Lemma is used to establish the following theorem.
\begin{theorem} {\bf (Convergence)}:\label{Theorem_con}
Suppose Assumptions \ref{cost-assump}--\ref{assum_strong_dual} hold, then for positive constant step-sizes satisfying:
 \eq{
\mu_w < {1 \over 2\delta-\nu} , \quad \mu_v < { \nu \over \lambda_{\max}(\cB\tran\cB)} \label{step-size-theorem}
}
recursions \eqref{primal-network_w} and \eqref{dual1-network}--\eqref{dual2-network} converge and it holds that $\sw_i$ converges to the optimal solution of \eqref{global-2}.
\end{theorem}
{\bf Proof}:  See Appendix \ref{Appendix:proof_theorem}. \qd

\noindent   At this point we showed that the dual coupled diffusion strategy, which handles non-smooth terms, converges to the optimal point.  However, it is still unclear how the sparsity of the constraints affects the convergence behavior. Apart from saving communication and memory, the next result reveals the advantage of exploiting the constraint structure. 
\begin{theorem} {\bf (Linear convergence)}:\label{coro-linear-conv}
Suppose Assumptions \ref{cost-assump}--\ref{assum_strong_dual} hold, and, furthermore, assume that each $R_k(w_k)=0$ and each matrix ${\rm blkcol}\{B_{e,k}\}_{e \in \cE_k}$ has full row rank. If the step sizes satisfy \eqref{step-size-theorem}, then  it holds that:
\eq{
&\|\tsw_i\|_{C_w}^2+\|\tsy_{i}\|^2_{\mu_w \over \mu_{v}}+ \|\tsx_{i}\|_{\mu_w \mu_v \Sigma}^2  \leq \gamma^i C_0
 \label{linear_corro}
}
for some constant $C_0>0$ where $C_w=I-\mu_w \mu_{v} \cB\tran \cB>0$ and
\eq{
\gamma \define \max \big\{\gamma_1,\gamma_2,1-\lambda_r \big\} < 1 }
with 
$\gamma_1=1-\mu_w \nu(1-\delta \mu_w)$, $\gamma_2=1-\mu_w \mu_v \lambda_{\min}(\cB \cB\tran)$ and  $\lambda_r$ denoting the smallest non-zero eigenvalue of $0.5(I-\cA)$.
\end{theorem}
{\bf Proof}:  See Appendix \ref{Appendix:proof_linear_conv}. \qd

\noindent The above result shows why solving \eqref{global-2} directly is important for at least two reasons. First, by using model \eqref{global-2}, we are able to prove linear convergence under the assumption that each ${\rm blkcol}\{B_{e,k}\}_{e \in \cE_k}$ has full row rank. If instead, we were to rewrite problem \eqref{global-2} into the form \eqref{global-1} by embedding zeros into the matrices $B_k$, then {\em our analysis} would require $B_k$ to be full row rank for linear convergence. This will not be satisfied if some agent is not involved in some constraint since in that case $B_k$ will have zero rows and, thus,  $B_k$ is row rank deficient even if ${\rm blkcol}\{B_{e,k}\}_{e \in \cE_k}$ has full row rank.

  The second more important reason is that the convergence rate depends on the connectivity of the sub-networks $\cC_e$ and not on the connectivity of the entire network, as we illustrate now. Note from the block structure of \eqref{calA} that the smallest non-negative eigenvalue of $0.5(I-\cA)$ has the form $
\lambda_r= \min_{e} \underline{\sigma}_e$
where $\underline{\sigma}_e$ denotes the smallest non-zero eigenvalue of the matrix $0.5(I-A_e)$. Since $
I-0.5(I-A_e)=0.5(I+A_e)=\bar{A_e}$, 
it holds that $
1-\underline{\sigma}_e= \bar{\lambda}_{2,e}$, where $\bar{\lambda}_{2,e}$ denotes the second largest eigenvalue of $\bar{A}_e$ (the largest eigenvalue is equal to one). Therefore,
\eq{
1-\lambda_r = 1- \min_{e} \underline{\sigma}_e =\max_{e} (1-\underline{\sigma}_e)= \max_{e} \bar{\lambda}_{2,e}
}  
Thus, assuming $1-\lambda_r$ is dominating the convergence rate, then the smaller $\max_{e} \bar{\lambda}_{2,e}$ is, the faster the algorithm is. We see that this depends on the second largest eigenvalue of the matrices $\{\bar{A}_e\}$, which depends on the sub-networks connectivity and not the whole network. This observation reveals the importance of the algorithm for sparse networks and under sparsely coupled constraints. Since in that case the small sub-networks are much well connected than the whole network. This observation will be illustrated in the simulation section next.
\begin{remark}[\sc Condition Number]{\em 
By using the  the  upper bound \eqref{step-size-theorem}, we conclude from Theorem \ref{coro-linear-conv} that the number of iterations needed to reach $\epsilon$ accuracy is on the order of
\eq{
O\left(\max\{\kappa_J \kappa_B, 1 / \lambda_r \} {\rm log} {1 \over \epsilon} \right)  \nonumber
}
 where $
\kappa_J=\delta/ \nu$ and  $\kappa_B=\lambda_{\max}(\cB \cB\tran)/\lambda_{\min}(\cB \cB\tran)$
are the condition numbers of the cost $\cJ(.)$ and the matrix $\cB\cB\tran$, respectively.
    }\qd
\end{remark}
%\begin{figure*}[!t]
%	\centering
%	\begin{subfigure}{.48\linewidth}
%		\includegraphics[scale=0.42]{ls_0917.pdf}
%		\caption{\small Least squares results.}
%		\label{fig:dule-coupled-diffusion-ls}
%	\end{subfigure} 
%	\hspace{5mm}
%	\begin{subfigure}{.48\linewidth}
%		\includegraphics[scale=0.42]{lr_0916.pdf}
%		\caption{\small Logistic regression results.}
%		\label{fig:dule-coupled-diffusion}
%	\end{subfigure} 
%	%\vspace{0.5mm}
%	\caption{ \small Simulation results. *Dual diffusion refers to \eqref{final-algorithm_jk=1} applied on the same problem reformulated into \eqref{global-1}, which ignores the sparsity structure. Similarly, both IDC-ADMM \cite{chang2015multi} and "dual DIGing" \cite{lee2017sublinear} are designed for problem \eqref{global-1} and ignore the sparsity structure.}
%	\label{fig:simulation-results}
%\end{figure*}
\begin{figure*}[!t]
	\centering
	\begin{subfigure}{.32\linewidth}
		\includegraphics[scale=0.33]{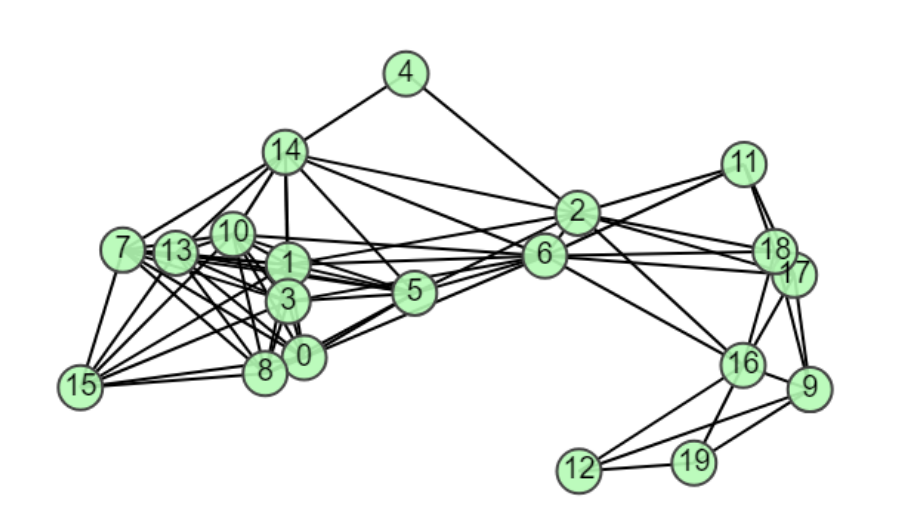}
		\vspace{8mm}
		\caption{\footnotesize Network topology used in simulations.} \label{fig:network}
	\end{subfigure}
	\begin{subfigure}{.3\linewidth}
		\includegraphics[scale=0.3]{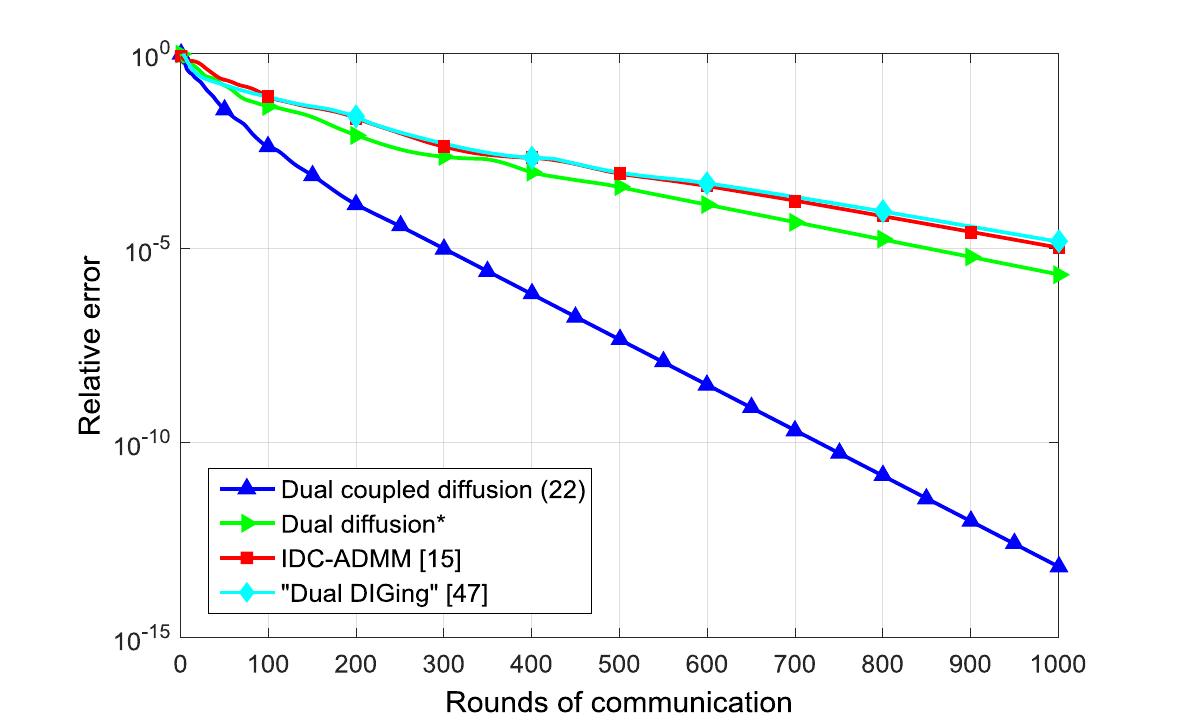}
		\caption{\footnotesize Least squares results.}
		\label{fig:dule-coupled-diffusion-ls}
	\end{subfigure} 
	\begin{subfigure}{.3\linewidth}
		\includegraphics[scale=0.3]{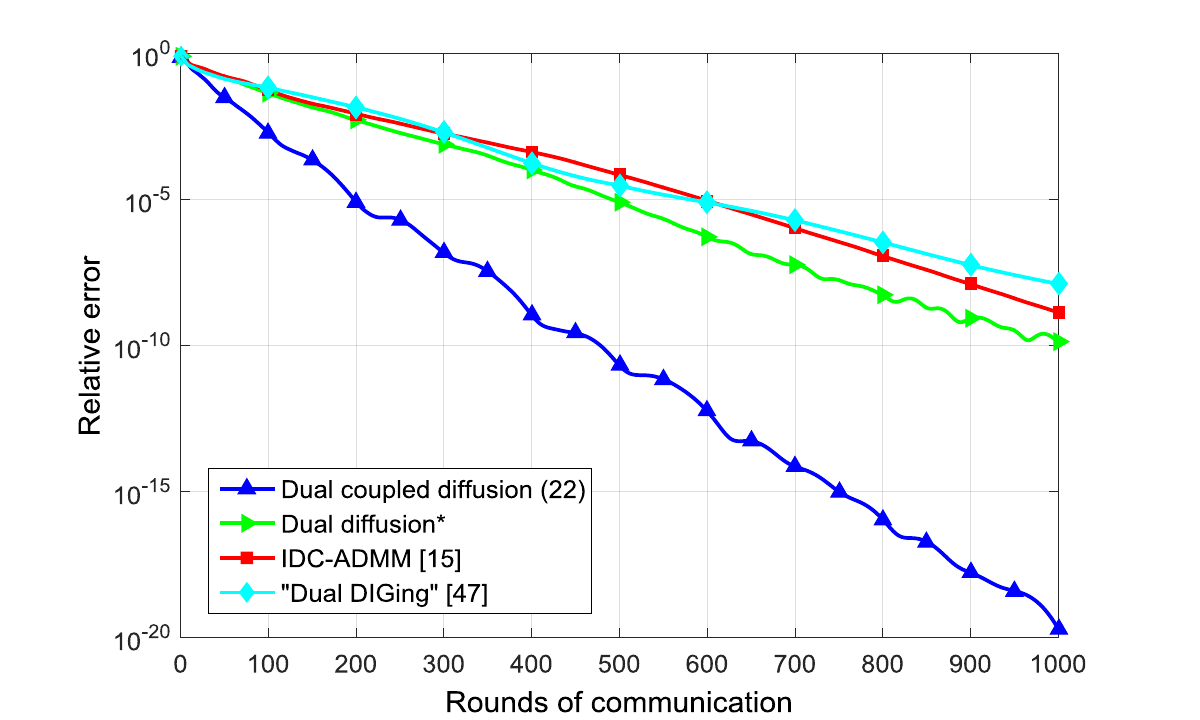}
		\caption{\footnotesize Logistic regression results.}
		\label{fig:dule-coupled-diffusion}
	\end{subfigure} 
	%\vspace{0.5mm}
	\caption{ \small Simulation results. *Dual diffusion refers to \eqref{final-algorithm_jk=1} applied on the same problem reformulated into \eqref{global-1}, which ignores the sparsity structure. Similarly, both IDC-ADMM \cite{chang2015multi} and "dual DIGing" \cite{lee2017sublinear} are designed for problem \eqref{global-1} and ignore the sparsity structure.}
	\label{fig:simulation-results}
\end{figure*}
\section{Numerical Simulation} \label{section:sim}
In this section, we test the performance of the proposed algorithm with two numerical experiments.
\begin{itemize}
\item {\em Distributed Linear Regression}:
The first set-up considers a linear regression problem with costs:
\eq{
\scalemath{0.95}{	J_k(w_k) = \frac{1}{2T_k}\sum_{t=1}^{T_k} \|u_{k,t}\tran w_k - p_{k}(t)\|^2} \nonumber 
}
and $R_k(w_k)=\eta_1 \|w_k\|_1$ where $u_{k,t} \in \real^{Q_k}$ is the regressor vector for data sample $t$, $p_k(t) \in \real$, and  $T_k$ denotes the amount of data for agent $k$.
\item {\em Distributed Logistic Regression}: The second set-up considers a logistic regression problem with costs:
\eq{
\scalemath{0.95}{ J_k(w_k) = \frac{1}{T_k}\sum_{t=1}^{T_k} \ln\left(1+\exp(-x_k(t) h_{k,t}\tran w_k)\right)+0.5\eta_2 \|w_k\|^2 } \nonumber
}
and $R_k(w_k)=\eta_1 \|w_k\|_1$. The vector $h_{k,t} \in \real^{Q_k}$ is the regressor vector for data sample $t$, and $x_k(t)$ is the label for that data sample, which is either $+1$ or $-1$.  
\end{itemize}
In both experiments, the network used is shown in Fig. \ref{fig:network} with $K=20$ agents. The positions ($x$-axis and $y$-axis) of the agents are randomly generated in $([0,1],[0,1])$,  and two agents are connected if the distance between them is less than or equal $d=0.3$. As for the constraints, we assume $E=K=20$, and each constraint $e$ (or $k$) (where $e\in \{1,\cdots,20\}$) is associated with a subnetwork involving agent $e$ (or $k$) and all its neighbors as described in equation \eqref{constraint-neighbor}. Each element in $B_{e,k}$ is generated according to the standard Gaussian distribution $\cN(0, 1)$. Each $b_{e,k}$ is also randomly generated and we guarantee that there exists a feasible solution to \eqref{global-2}. All the combination matrices are generated according to the Metropolis rule.
%\begin{figure}[H]
%\centering
%		\includegraphics[scale=0.5]{network_for_ls_lr_seed_2018.pdf}
%		\caption{\small The network topology used in simulations.}
%		\label{fig:network}
%\end{figure}	

In the first simulation, we set $T_k=1000$ for all $k$ and each regressor $u_{k,t}$ is generated according to the Gaussian distribution $\cN(0, 1)$. To generate the associated $p_k(t)$, we first generate a vector $w_{k,0}\in \real^{Q_k}$ randomly from $\cN(0, 1)$. We let $20\%$ of the entries of $w_{0,k}$ to be $0$. With such sparse $w_{k,0}$, we generate $p_k(t)$ as $p_k(t) = u_{k,t}\tran w_{k,0} + n_k$ where $n_k\sim \cN(0, 0.1)$ is some Gaussian noise. In this experiment, we set $Q_k=10$ for $k=1,\cdots, K$. We also set $\eta_1=0.3$ and $B_{e,k} \in \real^{3 \times 10}$ to be an under-determined coefficient matrix.  In the second set-up, each $T_k=1000$. Among all local data samples, half of them are generated by the Gaussian distribution $\cN(1, 1)$ and their corresponding labels $\{x_k(t)\}$ are $+1$'s. The other half are generated by $\cN(-1,1)$ and their corresponding labels $\{x_k(t)\}$ are $-1$'s. We set $Q_k=10$ for $k=1,\cdots, K$ and $\eta_1=\eta_2=0.1$.  We let $B_{e,k} \in \real^{3 \times 10}$ to be an under-determined coefficient matrix. 
\begin{figure*}[!t] 
	\centering 
		\includegraphics[scale=0.37]{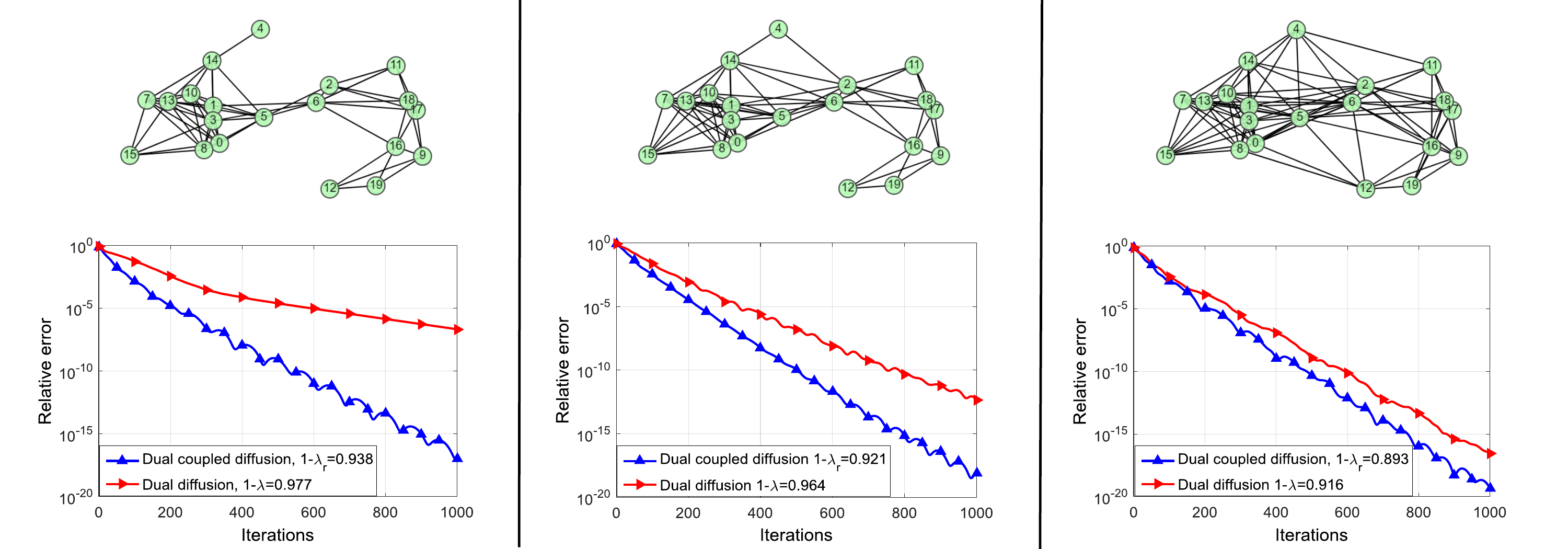}
		\caption{\small A comparison of algorithm \eqref{final-algorithm_jk=1} for different network connectivity and under two implementations: dual coupled diffusion exploits structure while dual diffusion ignores the structure.} \label{fig:sim_spars}
\end{figure*}
	
	To illustrate the effect of the constraint structure, we consider two approaches to solve problem \eqref{global-2}. The first
	approach is to use the dual coupled diffusion \eqref{final-algorithm_jk=1} while
	considering the structure of the problem \eqref{global-2}, i.e., run \eqref{final-algorithm_jk=1} with 
	$E = K, \cC_e = \cN_e$. The second approach is to ignore
	the special structure of the problem and reformulate it into
	the form of problem \eqref{global-1} and also run the dual coupled diffusion \eqref{final-algorithm_jk=1} with $E = 1, \cC_1 = \{1,\cdots,K\}$, which we call dual diffusion. To compare with other related methods that only share dual variables, we simulate the inexact distributed consensus ADMM (IDC-ADMM) from \cite{chang2015multi} and a modified proximal version of the one in \cite{lee2017sublinear} in which the dual iterates are updated similar to the DIGing algorithm in \cite{nedic2017achieving}, which we call ``Dual DIGing".   Both of these algorithm are designed for problem \eqref{global-1} and ignores any structure. The step-sizes are chosen manually to get the best possible performance for each algorithm. In the first linear regression setup, the parameters used are $(\mu_w=0.28, \ \mu_v=0.28)$ for the dual coupled diffusion, $(\mu_w=0.28, \ \mu_v=0.28)$ for the dual diffusion, $(c=0.25,\ \mu_w=0.05)$ for the IDC-ADMM  \cite{chang2015multi}, and the step-sizes are set to $0.45$ for the dual DIGing method. In the second logistic regression set-up, they are set to $(\mu_w=0.2, \ \mu_v=0.2)$ for the dual coupled diffusion, $(\mu_w=0.2, \ \mu_v=0.2)$ for the dual diffusion, $(c=0.45,\ \mu_w=0.2)$ for the IDC-ADMM  \cite{chang2015multi}, and the step-sizes are set to $0.18$ for the dual DIGing method. Figure \ref{fig:simulation-results} shows the relative error  $\frac{1}{K}\sum_{k=1}^K\|w_{k,i}-w_k^\star\|^2/\|w_k^\star\|^2$ for each of the previous algorithms for both set-ups. Note that the dual DIGing algorithm requires communicating two vectors each round of communication. It is observed that dual diffusion, the IDC-ADMM, and the dual DIGing algorithms have a close performance (all ignores any structure), while the dual coupled diffusion clearly outperforms them.  This means that, apart from requiring less amount of data to be exchanged per round of communication, our algorithm is also able to reach an $\epsilon$ accuracy (where $\epsilon$ is arbitrarily small) with much less time compared to these other algorithms.  As explained before, this superiority is due to the sub-networks being better connected compared to the whole network and the dual coupled diffusion takes advantage of that. In this simulation, we have $1-\lambda_r=0.911$ for the dual coupled diffusion and $1-\lambda=0.973$  for the dual diffusion (we dropped the sub-index since we have one network combination matrix in this case), which backs up our theoretical findings. 
	
	 To further illustrate the effect of the sub-networks connectivity on the convergence rate, we simulate the dual coupled diffusion (exploits sparsity) and dual diffusion (which does not exploit the sparsity) with the same logistic regression set-up from before but for the three different networks shown in top half of Fig. \ref{fig:sim_spars}. The step sizes used in this simulation are adjusted to get the best possible results, which are shown on the bottom of Figure \ref{fig:sim_spars}. Note that the network on the left has less connections compared to the network on the right, and thus, the sub-networks on the left are more sparse than the one on the right. Note further that for the constraints settings used \eqref{constraint-neighbor}, the more connections the network has, the closer the sub-networks are to the entire network. It is seen that  dual coupled diffusion performs significantly better under sparser networks since in that case the sub-networks are much better connected than the whole network. On the other hand, when we add more connections, the sub-networks connectivity  becomes closer to the network connectivity and, thus, the performance of the two algorithms become closer and closer. The performance will become identical when all agents are involved in all the constraint.
 \section{Concluding remarks}
This work developed a proximal diffusion strategy
with guaranteed exact convergence for a multi-agent optimization problem with multiple coupled constraints. We established
analytically, and by means of simulations,
the superior convergence properties of an algorithm that considers the sparsity structure in the constraints compared to others that ignore this structure.
\appendices
\section{Equivalent Representation}\label{Appendix:dual-represe}
In this appendix, we show that  \eqref{dual1-network}--\eqref{dual2-network} is equivalent to \eqref{networkrecursion2}. Multiplying equation \eqref{dual1-network} by $\cU_1 \Sigma$ and then collecting the term $\cU_1 \Sigma \ssx_{i-1}$ we get:
\eq{
\cU_1 \Sigma \ssx_{i} &=(I-\cU_1 \Sigma \cU_1\tran )\cU_1 \Sigma \ssx_{i-1}  \nonumber \\
& \quad - \mu_v^{-1}\left(\cU_1 \Sigma \cU_1\tran \right)  \bigg(\sy_{i-1}+\mu_v \big( \cB \sw_{i}-  b\big)  \bigg)
}
Let $\bar{\ssx}_i \define \cU_1 \Sigma \ssx_{i}$. Using \eqref{I-Adecomposition} and collecting the term $\ssx_{i-1}$ on the right hand side of the last equation, we get:
\eq{
\bar{\ssx}_{i}&=\bar{\cA} \bar{\ssx}_{i-1} -\mu_v^{-1}  \frac{1}{2} \left(I-\cA\right) \left(\sy_{i-1}+\mu_v \big( \cB \sw_{i}-  b\big)  \right) \label{appenid}}
Multiplying \eqref{dual2-network} by $\bar{\cA}$ on the left and using the definition $\bar{\ssx}_i \define \cU_1 \Sigma \ssx_{i}$ we have:
\eq{
\bar{\cA}  \sy_{i-1}&=\bar{\cA}  \sy_{i-2}+\mu_v \bar{\cA}  \big( \cB \sw_{i-1}-  b\big)+\mu_v \bar{\cA}  \bar{\ssx}_{i-1} \label{app1_1}
}
Now, subtracting \eqref{app1_1} from \eqref{dual2-network} we get:
\eq{
&\sy_{i} - \bar{\cA}  \sy_{i-1} \nonumber \\
&=\sy_{i-1} - \bar{\cA}  \sy_{i-2}+\mu_v  \big( \cB \sw_{i}-  b\big)-\mu_v \bar{\cA}  \big( \cB \sw_{i-1}-  b\big)\nonumber \\
& \ +\mu_v   (\bar{\ssx}_{i} - \bar{\cA} \bar{\ssx}_{i-1}) 
}
Using \eqref{appenid} we can remove the term $\mu_v   (\bar{\ssx}_{i} - \bar{\cA}  \bar{\ssx}_{i-1})$ from the previous expression to get: 
\eq{
&\sy_{i} - \bar{\cA}  \sy_{i-1} \nonumber \\
&=\sy_{i-1} - \bar{\cA}  \sy_{i-2}+\mu_v  \big( \cB \sw_{i}-  b\big)-\mu_v \bar{\cA}  \big( \cB \sw_{i-1}-  b\big)\nonumber \\
& \quad - \frac{1}{2} \left(I-\cA\right) \left(\sy_{i-1}+\mu_v \big( \cB \sw_{i}-  b\big)  \right) \nonumber \\
&=   \bar{\cA} \sy_{i-1} - \bar{\cA}  \sy_{i-2}+\mu_v  \bar{\cA} \big( \cB \sw_{i}-  b\big)-\mu_v \bar{\cA}  \big( \cB \sw_{i-1}-  b\big) \nonumber
}
Rearranging the last expression gives \eqref{networkrecursion2}.

\section{Primal Error Bound \eqref{werrrrr}} \label{Appendix:Primal_error}
 From the optimality condition of \eqref{primal-network_w}, we have:
\eq{
\sw_i&=\sw_{i-1}-\mu_w \grad\cJ(\sw_{i-1}) -\mu_w \cB\tran \sy_{i-1} -\mu_w g_i 
} 
for some $g_i \in \partial_{\ssw} \cR(\sw_{i})$. 
Rearranging the last equation and using the optimality condition \eqref{optimality1} we get:
\eq{
\sw_{i-1}-\sw_i&=\mu_w \big(\grad\cJ(\sw_{i-1})-\grad\cJ(\sw^\star)\big) +\mu_w (g_i-g^\star) \nonumber \\
& \quad +\mu_w \cB\tran (\sy_{i-1}-\sy^\star)
}
Multiplying $(\sw^\star-\sw_i)\tran$ to both sides of the previous equation, we get:
\eq{
&(\sw^\star-\sw_i)\tran(\sw_{i-1}-\sw_i) \nonumber \\&=\mu_w (\sw^\star-\sw_i)\tran \big(\grad\cJ(\sw_{i-1})-\grad\cJ(\sw^\star)\big) \nonumber \\
&\quad +\mu_w (\sw^\star-\sw_i)\tran(g_i-g^\star) \nonumber \\
& \quad +\mu_w (\sw^\star-\sw_i)\tran \cB\tran (\sy_{i-1}-\sy^\star)  \label{linearapp111}
}
From the conditions on $R_k(w_k)$ in Assumption \ref{cost-assump}, there exists at least one subgradient at every point. And from the subgradient property \eqref{subgradients}  we have $g_x\tran(y-x)\leq f(y)-f(x)$ and $g_y\tran(x-y)\leq f(x)-f(y)$. Summing the two inequalities with $y=\sw^\star$ and $x=\sw_i$, we get $(\sw^\star-\sw_i)\tran(g_i-g^\star) \leq 0$. Using this bound in \eqref{linearapp111} we get: 
\eq{
&(\sw^\star-\sw_i)\tran(\sw_{i-1}-\sw_i) \nonumber \\ &\leq \mu_w (\sw^\star-\sw_i)\tran \big(\grad\cJ(\sw_{i-1})-\grad\cJ(\sw^\star)\big) \nonumber \\
& \quad +\mu_w (\sw^\star-\sw_i)\tran \cB\tran (\sy_{i-1}-\sy^\star)  \label{linearapp1}
}
 Note that:
\eq{
&2(\sw^\star-\sw_i)\tran(\sw_{i-1}-\sw_i)\nonumber \\
&=\scalemath{0.95}{-\|\sw^\star-\sw_i-(\sw_{i-1}-\sw_i)\|^2+\|\sw^\star-\sw_i\|^2+\|\sw_{i-1}-\sw_i\|^2} \nonumber \\
&=-\|\tsw_{i-1}\|^2+\|\tsw_i\|^2+\|\sw_{i-1}-\sw_i\|^2
}
Substituting the last equation into \eqref{linearapp1} and rearranging terms gives:
\eq{
&\|\tsw_i\|^2 -\|\tsw_{i-1}\|^2 \nonumber \\
&\leq -\|\sw_{i-1}-\sw_i\|^2 -2\mu_w (\sy_{i-1}-\sy^\star)\tran \cB (\sw_i-\sw^\star) \nonumber \\
& \quad  -2\mu_w (\sw_i-\sw^\star)\tran \big(\grad\cJ(\sw_{i-1})-\grad\cJ(\sw^\star)\big)  \label{linearapp2}
}
Using Assumption \ref{cost-assump} we can bound the inner product:
\eq{
&(\sw_{i}-\sw^\star)\tran \grad\cJ(\sw_{i-1}) \nonumber \\
&= (\sw_{i}-\sw_{i-1}+\sw_{i-1}-\sw^\star)\tran \grad\cJ(\sw_{i-1}) \nonumber \\
&\overset{\eqref{stron-convexity}}{\geq} (\sw_{i}-\sw_{i-1})\tran \grad\cJ(\sw_{i-1}) \nonumber \\
& \quad +\cJ(\sw_{i-1})-\cJ(\sw^\star)+ {\nu \over 2} \|\tsw_{i-1}\|^2  \nonumber \\
&= (\sw_{i}-\sw_{i-1})\tran \big(\grad\cJ(\sw_{i-1})-\grad\cJ(\sw_{i})+\grad\cJ(\sw_{i})\big) \nonumber \\
& \quad +\cJ(\sw_{i-1})-\cJ(\sw^\star)+ {\nu \over 2}  \|\tsw_{i-1}\|^2  }
We again use \eqref{stron-convexity} in the last expression to get:
\eq{
&(\sw_{i}-\sw^\star)\tran \grad\cJ(\sw_{i-1}) \nonumber \\
&\geq (\sw_{i}-\sw_{i-1})\tran \big(\grad\cJ(\sw_{i-1})-\grad\cJ(\sw_{i})\big) \nonumber \\
& \ +\cJ(\sw_{i})-\cJ(\sw_{i-1})+ {\nu \over 2}  \|\sw_{i}-\sw_{i-1}\|^2 \nonumber \\
& \ +\cJ(\sw_{i-1})-\cJ(\sw^\star)+ {\nu \over 2}  \|\tsw_{i-1}\|^2   \nonumber \\
&= -(\sw_{i-1}-\sw_{i})\tran \big(\grad\cJ(\sw_{i-1})-\grad\cJ(\sw_{i})\big) \nonumber \\
& \ + {\nu \over 2}  \|\sw_{i}-\sw_{i-1}\|^2  +\cJ(\sw_{i})-\cJ(\sw^\star)+ {\nu \over 2}  \|\tsw_{i-1}\|^2  
}
From \eqref{stron-convexity} it holds that:
\eq{
(\sw_{i}-\sw^\star)\tran \grad \cJ(\sw^\star) \leq \cJ(\sw_{i})-\cJ(\sw^\star) -{\nu \over 2}  \|\tsw_i\|^2  }
Therefore, the last inner product in \eqref{linearapp2} can be bounded as follows:
\eq{
& -2\mu_w (\sw_i-\sw^\star)\tran \big(\grad\cJ(\sw_{i-1})-\grad\cJ(\sw^\star)\big) \nonumber \\
 &= -2\mu_w (\sw_i-\sw^\star)\tran \grad\cJ(\sw_{i-1}) +2\mu_w (\sw_i-\sw^\star)\tran\grad\cJ(\sw^\star) \nonumber \\
 & \leq 2\mu_w(\sw_{i-1}-\sw_{i})\tran \big(\grad\cJ(\sw_{i-1})-\grad\cJ(\sw_{i})\big) \nonumber \\
& \quad -\mu_w \nu \|\sw_{i}-\sw_{i-1}\|^2 -\mu_w \nu \|\tsw_{i-1}\|^2 -\mu_w \nu \|\tsw_i\|^2 \nonumber \\
 & \leq \mu_w(2 \delta-\nu)\|\sw_{i}-\sw_{i-1}\|^2   -\mu_w \nu (\|\tsw_{i-1}\|^2 + \|\tsw_i\|^2) \label{z1z1z1}
 }
 where the last step holds because $
(\sw-z)\tran \big(\grad \cJ(\sw)-\grad \cJ(z)\big)  \leq \delta \|\sw-z\| ^2 $
holds by using the Cauchy-Schwartz inequality and \eqref{lipschitz}. Substituting \eqref{z1z1z1} into \eqref{linearapp2} gives \eqref{werrrrr}.

\section{Dual Error Bound \eqref{error3a}} \label{Appendix:dual_error}
 It holds that:
\eq{
& \|\tsx_{i-1}\|^2_{\cD}+\|\tsy_{i-1}\|^2_{\mu_{v}^{-1}} \nonumber \\
=&\|\ssx^\star-\ssx_{i}+\ssx_{i}-\ssx_{i-1}\|_\cD^2+\|\sy^\star-\sy_{i}+\sy_{i}-\sy_{i-1}\|_{\mu_{v}^{-1}}^2 \nonumber \\
=&\|\ssx^\star-\ssx_{i}\|_\cD^2+\|\sy^\star-\sy_{i}\|_{\mu_{v}^{-1}}^2 + \|\ssx_{i}-\ssx_{i-1}\|_\cD^2 \nonumber \\
&+\|\sy_{i}-\sy_{i-1}\|_{\mu_{v}^{-1}}^2+
2(\ssx_{i}-\ssx_{i-1})\tran \cD(\ssx^\star-\ssx_{i})\nonumber \\
&+2(\sy_{i}-\sy_{i-1})\tran \mu_{v}^{-1}(\sy^\star-\sy_{i}) \label{error1}
}
Rearranging the last equality we have:
\eq{
&\|\tsy_{i}\|^2_{\mu_{v}^{-1}}+\|\tsx_{i}\|_\cD^2-\|\tsy_{i-1}\|_{\mu_{v}^{-1}}^2-\|\tsx_{i-1}\|_\cD^2   \nonumber \\
&=  -\|\ssx_{i}-\ssx_{i-1}\|_\cD^2-\|\sy_{i}-\sy_{i-1}\|^2_{\mu_{v}^{-1}} \nonumber \\
& \ \ + 
2(\ssx_{i-1}-\ssx_i)\tran \cD(\ssx^\star-\ssx_{i})-2(\sy_{i-1}-\sy_i)\tran {\mu_{v}^{-1}}(\sy_{i}-\sy^\star) \label{error2}
}
Note that:
\eq{
&(\sy_{i}-\sy^\star)\tran \cU_1 \Sigma (\ssx_i-\ssx^\star) \nonumber \\
&\overset{(a)}{=} (\sy_{i}-\sy_{i-1}+\sy_{i-1}-\sy^\star)\tran \cU_1 \Sigma(\ssx_i-\ssx^\star)\nonumber \\
&\overset{(b)}{=}(\cU_1\tran(\sy_{i}-\sy_{i-1})+\cU_1\tran\sy_{i-1})\tran  \Sigma(\ssx_i-\ssx^\star)\nonumber \\
&\overset{\eqref{dual2-network}}{=} \big(\cU_1\tran(\mu_{v} (\cB \sw_i-b) +\mu_{v} \cU_1 \Sigma \ssx_i) +\cU_1\tran\sy_{i-1}\big)\tran \Sigma(\ssx_i-\ssx^\star) \nonumber
\\
&= \bigg( \cU_1\tran \big(\sy_{i-1}+\mu_{v}(\cB \sw_i-b)+\mu_{v} \cU_1 \Sigma \ssx_{i-1} \big) \nonumber \\
& \quad +\mu_{v}  \Sigma(\ssx_i-\ssx_{i-1}) \bigg)\tran  \Sigma(\ssx_i-\ssx^\star) \nonumber \\
&\overset{\eqref{dual1-network}}{=} \bigg(\mu_v (\ssx_{i-1}-\ssx_i) +\mu_{v}\Sigma(\ssx_{i}-\ssx_{i-1})\bigg)\tran  \Sigma(\ssx_i-\ssx^\star)
\nonumber \\
&= -\big(\mu_v( I-\Sigma )(\ssx_i-\ssx_{i-1})\big)\tran \Sigma (\ssx_i-\ssx^\star)
\nonumber \\
&\overset{\eqref{cal_D}}{=} -(\ssx_i-\ssx_{i-1})\tran \cD (\ssx_i-\ssx^\star) 
}
where in step (b) we took $\cU_1$ inside the first bracket and used $\cU_1\tran \sy^{\star}=0$ from \eqref{optimality3}.  From  step (a) and the last step we get:
\eq{
&(\sy_{i-1}-\sy^\star)\tran \cU_1 \Sigma(\ssx_i - \ssx^\star) \nonumber \\
&=
 -(\sy_i-\sy_{i-1})\tran  \cU_1 \Sigma(\ssx_i-\ssx^\star) -(\ssx_i-\ssx_{i-1})\tran \cD (\ssx_i-\ssx^\star) 
\label{aaa}}
Furthermore, note that:
\eq{
&(\sy_{i-1}-\sy^\star)\tran (\cB \sw_{i}-  b-\cB \sw^\star+  b)  \nonumber \\
&\overset{\eqref{dual2-network}}{=} (\sy_{i-1}-\sy^\star)\tran \left(- \mu_{v}^{-1}(\sy_{i-1}-\sy_i)-\cU_1 \Sigma \ssx_i-\cB \sw^\star+  b\right) \nonumber \\
&\overset{\eqref{optimality2}}{=} (\sy_{i-1}-\sy^\star)\tran \left(- \mu_{v}^{-1}(\sy_{i-1}-\sy_i)-\cU_1 \Sigma \ssx_i+\cU_1 \Sigma\ssx^\star\right) \nonumber \\
&= -(\sy_{i-1}-\sy^\star)\tran  \mu_{v}^{-1}(\sy_{i-1}-\sy_i)\nonumber \\
& \quad -(\sy_{i-1}-\sy^\star)\tran \cU_1 \Sigma(\ssx_i-\ssx^\star) \label{aaa1}
} 
Substituting \eqref{aaa} into \eqref{aaa1}, we have
\eq{
&(\sy_{i-1}-\sy^\star)\tran (\cB \sw_{i}-  b-\cB \sw^\star+  b)  \nonumber \\
&= -(\sy_{i-1}-\sy^\star)\tran \mu_{v}^{-1}(\sy_{i-1}-\sy_i) \nonumber \\
& \quad +(\sy_i-\sy_{i-1})\tran  \cU_1 \Sigma(\ssx_i-\ssx^\star)+(\ssx_i-\ssx_{i-1})\tran \cD (\ssx_i-\ssx^\star) 
 \nonumber \\
&= \left(-\mu_{v}^{-1}(\sy_{i-1}-\sy^\star)-\cU_1 \Sigma (\ssx_i-\ssx^\star)\right)\tran (\sy_{i-1}-\sy_i) \nonumber \\
& \quad +(\ssx_i-\ssx_{i-1})\tran \cD (\ssx_i-\ssx^\star) \nonumber \\
&\overset{(a)}{=}\bigg(-\mu_{v}^{-1}(\sy_{i-1}-\sy^\star)+\mu_{v}^{-1}(\sy_{i-1}-\sy_i) \nonumber \\
& \quad +(\cB \sw_{i}-  b)-\cB \sw^\star+  b\bigg)\tran (\sy_{i-1}-\sy_i)\nonumber \\
&\quad \quad+(\ssx_i-\ssx_{i-1})\tran \cD (\ssx_i-\ssx^\star) \nonumber \\
&= \left(-\mu_{v}^{-1}(\sy_{i}-\sy^\star)+\cB (\sw_i-\sw^\star)\right)\tran (\sy_{i-1}-\sy_i)\nonumber \\
& \quad +(\ssx_i-\ssx_{i-1})\tran \cD (\ssx_i-\ssx^\star) \nonumber \\
&= -(\sy_{i}-\sy^\star)\tran \mu_{v}^{-1}(\sy_{i-1}-\sy_i) + (\sw_i-\sw^\star)\tran \cB\tran (\sy_{i-1}-\sy_i)\nonumber \\
& \quad +(\ssx_i-\ssx_{i-1})\tran \cD (\ssx_i-\ssx^\star) \label{bbb}
} 
where in step (a) we used \eqref{dual2-network} and the optimality condition \eqref{optimality2}. Re-arranging the last equation \eqref{bbb}, we get
\eq{
& -(\sy_{i-1}-\sy_i)\tran \mu_v^{-1}(\sy_{i}-\sy^\star) +(\ssx_{i-1}-\ssx_i)\tran \cD (\ssx^\star-\ssx_i) \nonumber \\
&= (\sy_{i-1}-\sy^\star)\tran \cB (\sw_{i}- \sw^\star)-(\sy_{i-1}-\sy_i)\tran \cB (\sw_i-\sw^\star)  \nonumber
}
Substituting the previous equation into \eqref{error2}, we get,
\eq{
&\|\tsy_{i}\|^2_{\mu_{v}^{-1}}+\|\tsx_{i}\|_\cD^2-\|\tsy_{i-1}\|_{\mu_{v}^{-1}}^2-\|\tsx_{i-1}\|_\cD^2   \nonumber \\
&=  -\|\ssx_{i}-\ssx_{i-1}\|_\cD^2-\|\sy_{i}-\sy_{i-1}\|^2_{\mu_{v}^{-1}}  \nonumber \\
& \ + 2(\sy_{i-1}-\sy^\star)\tran \cB (\sw_{i}- \sw^\star)   -2(\sy_{i-1}-\sy_i)\tran \cB (\sw_i-\sw^\star)  \label{error3}}
The last term of \eqref{error3} can be rewritten as:
\eq{
&-2(\sy_{i-1}-\sy_i)\tran \cB (\sw_i-\sw^\star)\nonumber \\
 &=-\|\sy_{i-1}-\sy_i+\mu_{v} \cB (\sw_i-\sw^\star) \|_{\mu_{v}^{-1}}^2 \nonumber \\
& \quad +\|\sy_{i-1}-\sy_i\|_{\mu_{v}^{-1}}^2+ \|\cB (\sw_i-\sw^\star)\|_{\mu_{v}}^2 \nonumber \\
&= -\|\Sigma( \ssx^\star- \ssx_i) \|_{\mu_{v}}^2   +\|\sy_{i-1}-\sy_i\|_{\mu_{v}^{-1}}^2+ \|\cB (\sw_i-\sw^\star)\|_{\mu_{v}}^2 \nonumber
}
where in the last step we used \eqref{dual2-network}, \eqref{optimality2}, and $\cU_1\tran \cU_1=I$. Substituting the last equality into \eqref{error3}, we get \eqref{error3a}.
\section{Proof of Theorem \ref{Theorem_con}} \label{Appendix:proof_theorem}
Let us introduce the quantity:
\eq{
V(\tsw_i,\tsy_{i},\tsx_{i})=\|\tsw_i\|^2+\mu_w  \big(\|\tsy_{i}\|^2_{\mu_{v}^{-1}}+\|\tsx_{i}\|_\cD^2 \big)
} 
Using \eqref{werrrrr}--\eqref{error3a} and $
\|\cB \tsw_i\|_{\mu_{v}}^2 \leq \mu_{v} \lambda_{\max}(\cB\tran \cB) \| \tsw_i\|^2$,
it holds that:
\eq{
&V(\tsw_{i},\tsy_{i},\tsx_{i})-V(\tsw_{i-1},\tsy_{i-1},\tsx_{i-1}) \nonumber \\
& \leq -\underbrace{(1+\mu_w \nu-2\mu_w \delta)}_{>0}\|\sw_i-\sw_{i-1}\|^2 \nonumber \\
& \quad - \mu_w \nu \|\tsw_{i-1}\|^2 - \mu_w \underbrace{\left( \nu-\mu_{v} \lambda_{\max}(\cB\tran \cB)\right)}_{>0} \|\tsw_i\|^2  \nonumber \\
& \quad  -\mu_w\|\ssx_{i}-\ssx_{i-1}\|_\cD^2-\mu_w\| \Sigma\tsx_i  \|_{\mu_{v}}^2  \leq 0 \label{diff-lyap}
}
where the last inequality holds under \eqref{step-size-theorem}. Since $V(\tsw_{i},\tsy_{i},\tsx_{i})$ is non-negative, we conclude that the norm of the error is non-increasing and bounded. Iterating the above inequality we have:
\eq{
&V(\tsw_{i},\tsy_{i},\tsx_{i}) \leq V(\tsw_{0},\tsy_{0},\tsx_{0})  \nnb
& - \sum_{j=1}^{i} \bigg((1+\mu_w \nu-2\mu_w \delta)\|\sw_j-\sw_{j-1}\|^2  + \mu_w \nu \|\tsw_{j-1}\|^2  \nonumber \\
& \quad+ \mu_w \left( \nu-\mu_{v} \lambda_{\max}(\cB\tran \cB)\right) \|\tsw_j\|^2 +\mu_w\|\ssx_{j}-\ssx_{j-1}\|_\cD^2 \nonumber \\
& \quad  +\mu_w\| \Sigma\tsx_j \|_{\mu_{v}}^2 \bigg)
}
and thus
\eq{
&\sum_{i=1}^{\infty} \bigg((1+\mu_w \nu-2\mu_w \delta)\|\sw_i-\sw_{i-1}\|^2  + \mu_w \nu \|\tsw_{i-1}\|^2  \nonumber \\
& \quad+ \mu_w \left( \nu-\mu_{v} \lambda_{\max}(\cB\tran \cB)\right) \|\tsw_i\|^2 +\mu_w\|\ssx_{i}-\ssx_{i-1}\|_\cD^2 \nonumber \\
& \quad  +\mu_w\| \Sigma\tsx_i \|_{\mu_{v}}^2 \bigg) 
\leq V(\tsw_{0},\tsy_{0},\tsx_{0})
}
Since the sum of the infinite positive terms is upper bounded by a constant, it holds that each term  $(\sw_i-\sw_{i-1}),\tsw_{i-1},\tsw_i,(\ssx_{i}-\ssx_{i-1})$, and $ \Sigma\tsx_i  $ must converge to zero.
\section{Proof of Theorem \ref{coro-linear-conv}} \label{Appendix:proof_linear_conv}
\begin{figure}[H]
	\centering
	\includegraphics[scale=0.4]{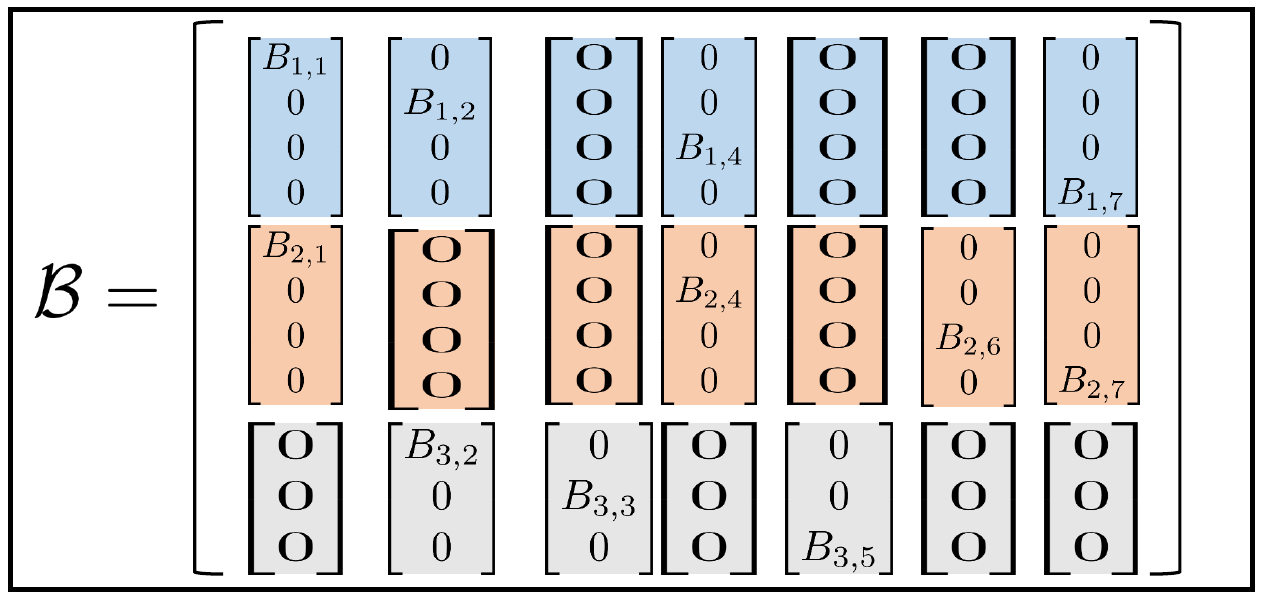}
	\caption{\footnotesize An illustration of the construction $\cB$ for the network in Figure \ref{fig:network-dual2}.}
	\label{fig:constructionBexample}
\end{figure}

 From the structure of $\cB$ in \eqref{calBtran}, it can be confirmed that $\cB$ having full row rank is equivalent to assuming that each matrix ${\rm blkcol}\{B_{e,k}\}_{e \in \cE_k}$ has full row rank. This is illustrated in Fig. \ref{fig:constructionBexample}. Because two different agents belonging to the same cluster are located differently in $\sy^e$, it holds that the block rows of $\cB$ are zeros except at one location. Recall that $B_{ek} \in \real^{S_e \times Q_k}$. Therefore, an equivalent statement is to say that  ${\rm blkcol}\{B_{e,k}\}_{e \in \cE_k}$ has full row rank. 
The last term in \eqref{linearapp2} can be rewritten as
\eq{
&-2\mu_w (\sw_i-\sw^\star)\tran \big(\grad\cJ(\sw_{i-1})-\grad\cJ(\sw^\star)\big) \nnb
& = -2\mu_w (\sw_{i-1}-\sw^\star)\tran \big(\grad\cJ(\sw_{i-1})-\grad\cJ(\sw^\star)\big) \nnb
& \quad +2\mu_w (\sw_{i-1}-\sw_i)\tran \big(\grad\cJ(\sw_{i-1})-\grad\cJ(\sw^\star)\big)  
}
the last term can be upper bounded by
\eq{
&2\mu_w (\sw_{i-1}-\sw_i)\tran \big(\grad\cJ(\sw_{i-1})-\grad\cJ(\sw^\star)\big) \nnb
& =  - \|\sw_{i-1}-\sw_i-\mu_w \big(\grad\cJ(\sw_{i-1})-\grad\cJ(\sw^\star)\big)\|^2 \nnb
& \quad +\|\sw_{i-1}-\sw_i\|^2 +\mu_w^2\|\grad\cJ(\sw_{i-1})-\grad\cJ(\sw^\star)\|^2 \nnb
& \overset{(a)}{=} - \mu_w^2 \|\cB\tran \tsy_{i-1}\|^2 + \|\sw_{i-1}-\sw_i\|^2  \nnb
& \quad +\mu_w^2\|\grad\cJ(\sw_{i-1})-\grad\cJ(\sw^\star)\|^2 \nnb
& \leq  - \mu_w^2 \|\cB\tran \tsy_{i-1}\|^2 + \|\sw_{i-1}-\sw_i\|^2  \nnb
& \quad +\mu_w^2 \delta
  (\sw_{i-1}-\sw^\star)\tran \big(\grad\cJ(\sw_{i-1})-\grad\cJ(\sw^\star)\big)
}
where in step (a) we used \eqref{primal-network_w} and \eqref{optimality1} with $\cR(\sw)=0$. The last inequality holds from \cite[Theorem 2.1.5]{nesterov2013introductory} since $\cJ(\sw)$ has $\delta$-Lipschitz gradients. Combining the last two equations we have
\eq{
&-2\mu_w (\sw_i-\sw^\star)\tran \big(\grad\cJ(\sw_{i-1})-\grad\cJ(\sw^\star)\big) \nnb
&\leq  - \mu_w^2 \|\cB\tran \tsy_{i-1}\|^2 + \|\sw_{i-1}-\sw_i\|^2  \nnb
& \ -\mu_w (2-\delta \mu_w)
  (\sw_{i-1}-\sw^\star)\tran \big(\grad\cJ(\sw_{i-1})-\grad\cJ(\sw^\star)\big) \nnb
&\leq  - \mu_w^2 \|\cB\tran \tsy_{i-1}\|^2 + \|\sw_{i-1}-\sw_i\|^2   -\mu_w \nu(2-\delta \mu_w)
  \|\tsw_{i-1}\|^2 \nonumber
} 
where the last step holds from the strong-convexity condition \eqref{stron-convexity} and $(2-\delta \mu_w)>0$ for $\mu_w < 2/\delta$. 
Substituting into \eqref{linearapp2} we get: 
\eq{
\|\tsw_i\|^2  &\leq  \big(1-\mu_w \nu(2-\delta \mu_w)\big) \|\tsw_{i-1}\|^2 - \mu_w^2 \|\cB\tran \tsy_{i-1}\|^2  \nnb
& \  -2\mu_w (\sy_{i-1}-\sy^\star)\tran \cB (\sw_i-\sw^\star)  
 \label{primal_linear_bound}
} 
Note that $-\mu_w\|\ssx_{i}-\ssx_{i-1}\|_\cD^2 \leq 0$. Thus, multiplying \eqref{error3a} by $\mu_w$ and rearranging terms we get:
\eq{
&\mu_w\|\tsy_{i}\|^2_{\mu_{v}^{-1}}+\mu_w \|\tsx_{i}\|_{\mu_v \Sigma}^2 =\mu_w\|\tsy_{i}\|^2_{\mu_{v}^{-1}}+\mu_w \|\tsx_{i}\|_{\cD+\mu_v \Sigma^2}^2  \nonumber \\
&\leq \mu_w \|\cB \tsw_i\|_{\mu_{v}}^2+ \mu_w\|\tsy_{i-1}\|_{\mu_{v}^{-1}}^2+\mu_w \|\tsx_{i-1}\|_\cD^2     \nonumber \\
&\quad + 2\mu_w(\sy_{i-1}-\sy^\star)\tran \cB (\sw_{i}- \sw^\star) 
\label{dual_linear_bound}
}
Since $\cB$ is full row rank, it holds that $
\|\cB\tran \tsy_{i-1}\|^2 \geq \lambda_{\min}(\cB \cB\tran)  \| \tsy_{i-1}\|^2$. Using this bound,  $\cD=\mu_v (\Sigma-\Sigma^2)$, and combining \eqref{primal_linear_bound} and \eqref{dual_linear_bound}, we get:
\eq{
&\|\tsw_i\|^2+\|\tsy_{i}\|^2_{\mu_w \over \mu_{v}}+ \|\tsx_{i}\|_{\mu_w \mu_v \Sigma}^2  \leq  \|\tsx_{i-1}\|^2_{\mu_w \mu_v (\Sigma-\Sigma^{2})}  \nonumber \\
& \ +  \big(1-\mu_w \nu(2-\delta \mu_w)\big) \|\tsw_{i-1}\|^2   
 + \mu_w \mu_{v} \|\cB \tsw_i\|^2 \nnb
 & \ + \big(1-\mu_w \mu_v \lambda_{\min}(\cB \cB\tran)\big)\|\tsy_{i-1}\|_{\mu_w \over \mu_{v}}^2  
 \label{linearbounda}
}
Since $\Sigma>0$ we have $
-\|\tsx_{i-1}\|^2_{\mu_w \mu_v \Sigma^{2}} \leq -\lambda_r \|\tsx_{i-1}\|^2_{\mu_w \mu_v \Sigma 
}$. 
Substituting this bound into \eqref{linearbounda} and rearranging, we arrive at the following inequality: 
\eq{
&\|\tsw_i\|_{C_w}^2+\|\tsy_{i}\|^2_{\mu_w \over \mu_{v}}+ \|\tsx_{i}\|_{\mu_w \mu_v \Sigma}^2 \leq  (1-\lambda_r)\|\tsx_{i-1}\|^2_{\mu_w \mu_v \Sigma}   \nonumber \\
& \ +  \big(1-\mu_w \nu(2-\delta \mu_w)\big) \|\tsw_{i-1}\|^2  +\gamma_2 \|\tsy_{i-1}\|_{\mu_w \over \mu_{v}}^2  
 \label{linearboundb}
}
where $C_w=I-\mu_w \mu_{v} \cB\tran \cB$ and $\gamma_2=1-\mu_w \mu_v \lambda_{\min}(\cB \cB\tran)$. Let $\gamma_1=\big(1-\mu_w \nu(1-\delta \mu_w)\big)$ and note that
\eq{
&\big(1-\mu_w \nu(2-\delta \mu_w)\big) \|\tsw_{i-1}\|^2 \nnb
&=\gamma_1 \|\tsw_{i-1}\|_{C_w}^2 -\mu_w \nu  \|\tsw_{i-1}\|^2  + \gamma_1 \mu_w \mu_{v}  \| \tsw_{i-1}\|_{\cB\tran \cB}^2 \nnb
&\leq\gamma_1 \|\tsw_{i-1}\|_{C_w}^2 -\mu_w \big(\nu-\mu_v\lambda_{\max}(\cB\tran \cB)\big) \|\tsw_{i-1}\|^2 \nnb
& \leq \gamma_1 \|\tsw_{i-1}\|_{C_w}^2 \nonumber
}
where the first inequality holds since $\gamma_1<1$  for $\mu_w < {1 \over 2\delta -\nu}  \leq {1 \over \delta}$ and the last inequality holds under \eqref{step-size-theorem}. Substituting into \eqref{linearboundb} we get:
\eq{
&\|\tsw_i\|_{C_w}^2+\|\tsy_{i}\|^2_{\mu_w \over \mu_{v}}+ \|\tsx_{i}\|_{\mu_w \mu_v \Sigma}^2 \nnb
& \leq  \gamma_1 \|\tsw_{i-1}\|_{C_w}^2  +\gamma_2 \|\tsy_{i-1}\|_{\mu_w \over \mu_{v}}^2  + (1-\lambda_r)\|\tsx_{i-1}\|^2_{\mu_w \mu_v \Sigma} \nonumber
}
Under condition \eqref{step-size-theorem}, it holds that $\mu_w \mu_v < 1/ \lambda_{\max}(\cB\tran\cB)$; thus, $\gamma_2=1-\mu_w \mu_v \lambda_{\min}(\cB \cB\tran)<1$ and $C_w=I-\mu_w \mu_{v} \cB\tran \cB>0$.   Since $0<\lambda_r<1$, we have $1-\lambda_r<1$.  By iterating the previous inequality   we arrive at \eqref{linear_corro}.

{
\bibliographystyle{IEEEtran}
\bibliography{IEEEabrv,myref_review}
}

\end{document}